**Evolution of Leibniz's thought in the matter of fictions and infinitesimals**

**Monica Ugaglia**
**Il Gallo Silvestre, Località Collina 38, Montecassiano, Italy**
**email: monica.ugaglia@gmail.com**

**Mikhail Katz**
**Department of Mathematics, Bar Ilan University, Ramat Gan 5290002, Israel**
**email: katzmik@math.biu.ac.il**

*Abstract.*
In this paper we offer a reconstruction of the evolution of Leibniz's thought concerning the problem of the infinite divisibility of bodies, the tension between actuality, unassignability and syncategorematicity, and the closely related question of the possibility of infinitesimal quantities, both in physics and in mathematics.

   Some scholars have argued that syncategorematicity is a mature acquisition, to which Leibniz resorts to solve the question of his infinitesimals – namely the idea that infinitesimals are just signs for Archimedean exhaustions, and their unassignability is a nominalist maneuver. On the contrary, we show that sycategorematicity, as a traditional idea of classical scholasticism, is a feature of young Leibniz's thinking, from which he moves away in order to solve the same problem, as he gains mathematical knowledge.

   We have divided Leibniz's path toward his mature view of infinitesimals into five phases, which are especially significant for reconstructing the entire evolution. In our reconstruction, an important role is played by Leibniz's text *De Quadratura Arithmetica*. Based on this and other texts we dispute the thesis that fictionality coincides with syncategorematicity,[1] and that unassignability can be bypassed. On the contrary, we maintain that unassignability, as incompatible with the principle of harmony, is the ultimate reason for the fictionality of infinitesimals.

1. INTRODUCTION

Leibniz's thinking about infinitely small quantities underwent a complex evolution, both from the mathematical and from the epistemological point of view.

   In his maturity, Leibniz was aware of the difficulties that the notion of an infinitely small quantity raises, but his writings suggest that his epistemological awareness has matured as a result of his mathematical proficiency. While the young Leibniz applied almost naively in physics the new mathematical discoveries of his time, of which he had only a superficial knowledge, in his maturity he develops a more complex view of the relationship between mathematics and physics (and metaphysics), and becomes more cautious both in the application and in the definition of mathematical notions.

   In particular, Leibniz eventually comes to master the notion of unassignability, and the consequences it involves. This leads him to a cogent treatment of infinitesimals in mathematics – infinitesimals are non-contradictory notions – and to a rejection in physics – they are incompatible with the principle of harmony. In other words, this leads Leibniz to his mature view of infinitesimals as accidental impossibilities, or "well-founded fictions."

---

[1]   In this paper we employ "syncategorematic" as a shorthand for "eventually identifiable with a procedure of exhaustion" and, as a consequence, "involving only assignable quantities". We also identify syncategorematicity with potentiality, as suggested by some part of the Scholastics, and in section 2.2 we show that the two characterisations are equivalent, provided that "potentiality" is intended in the correct way: namely, as in the unending iterative procedures of Greek mathematics.

We have divided Leibniz's path toward this position into five phases, which can be roughly summarized as follows.

The young Leibniz shares with many contemporary philosophers both a scholastic background, and the need of correcting it in the wake of the modern mechanistic philosophy. At this stage, Leibniz's thought does not require infinitely small quantities, for both in mathematics and in physics he adopts a traditional syncategorematic approach: both extension and matter are continuous, which means that they are potentially infinitely divisible, and every obtained part is assignable. Take a continuous quantity – for instance a segment – you can divide it, and then divide it, and then divide it... and so on, indefinitely. Since the procedure has no end, no ultimate constituents can be reached, in either mathematics or physics, and since the quantity obtained at every step of the procedure is smaller than the one obtained at the previous step, but is nonetheless finite, no unassignable infinitesimal quantities can be reached, either (Phase I: 1668-1669).

Quite early, for metaphysical reasons and in opposition to Descartes, the question of the ultimate constituents of physical bodies becomes a matter of primary importance for Leibniz. However, the solution he offers before his stay in Paris – marking the start of Leibniz's mathematical training and his acquaintance with the contemporary debate – is a very problematic one. At the very beginning (phase II: 1670-1671) Leibniz simply redefines matter as a continuous stuff infinitely divided – and not only divisible – into parts. Both actuality and unity of these parts are due to their *conatus*, or unextended amounts of motion. The geometrical counterpart of this physical situation is a hybrid in which the infinite divisibility of the continuum has to coexist with the actuality of unextended "fat points", which Leibniz identifies with Cavalieri's indivisibles. In this phase unassignability is accepted, without any further elaboration.

In the first years of the Parisian stay (phase III: 1672-1674), a slightly modified version of this solution is devised, which seems to solve the problem of actuality in physics, but still raises difficulties in mathematics.[2] The core idea is to put discreteness before continuity: physical bodies do not come from an infinitely divided continuous stuff, but from the aggregation of discrete basic elements. As in phase II the ultimate elements are actually infinite and "glued" together by their *conatus* of motion. In mathematics too the situation is similar, with the only difference that now points are extended, although unassignable. Apparently, in this phase to keep together extended but unassignable quantities and assignable ones is not a problem for Leibniz. Unassignable quantities are not called fictions, nor referred to using cognate verbal forms.

Only towards the end of his stay in Paris, Leibniz becomes aware of the existence of a problem, as the notes he wrote in phase IV (1675-1676) show, as well as his further writings. In particular, Leibniz rejects the possibility that physical quantities are not all comparable to each other: a world containing some portion of matter which has an unassignable ratio to another portion would be less perfect than a world where all portions are comparable.

For this reason, the solution Leibniz offers in this period to the problem of the composition of physical bodies is only apparently similar to the one described in phase III. In both cases we have an infinitely divided material stuff, and a sort of glue that keeps the parts together, but while in phase III the glue, namely the *conatus* of motion, corresponded to physical objects, in phase IV material parts are held together by metaphysical bonds. In other words, in phase III a geometrical description was needed, which did contain both the parts of a body, which are assignable, although infinitely divisible, and their unassignable *conatus*. In phase IV, instead, a geometrical description is needed only for the assignable physical parts, while the metaphysical "cement" is no longer constrained.

---

2   See further in note 21.

It is exactly in phase IV that Leibniz starts referring to unassignable quantities as "fictions" and links them to other "imaginary" notions, such as imaginary numbers, that are non-contradictory but do not admit any traditional geometrical representation.

Finally (phase V: 1677- ...) the separation between mathematics, physics and metaphysics is gradually brought to completion. In particular, mathematics ceases to be the foundation of physics, and geometry ceases to be its natural description: the true essence of physics is metaphysics, and as such it is non representable. In particular, being extraneous to the category of quantity, metaphysical objects can be thought of as unassignable in a stronger sense, and without problematic consequences.

## 2. The quest for actuality

### 2.1 Actuality and the ultimate elements of physics

In the very first years of his philosophical career, here described as phase I, Leibniz still operates in a scholastic context, where physical objects are continuous and continuity implies potentiality, by definition. But as soon as he begins to develop his own positions, the need to go beyond potentiality, grounding physical objects on something actual, becomes a key feature of Leibniz's philosophy.

The search for actuality, however, is forced to coexist with another crucial feature, namely the infinite divisibility of matter, which Leibniz inherits from scholasticism and will never abandon. The tension between these two seemingly irreconcilable positions, combined with a firm belief in keeping metaphysics out of physics, is what characterizes Leibniz's first moves towards an original theory of physical substance.

In the early 1670s, here described as phases II and III, we witness Leibniz's first unsuccessful attempts to conciliate actuality and infinite divisibility, without resorting to metaphysics. The idea was to locate actuality in some ultimate elements of matter, making it derive not from a vague metaphysical form, but from a much more concrete, and mathematically manageable, (tendency to) motion. In Leibniz's intentions, by virtue of their motion, the ultimate parts of a physical body should be able to stay connected together, as a whole, while maintaining their individuality. But neither the continuous parts of phase II, nor the discrete atoms of phase III were up to the task.

At this point Leibniz shifts his attention from material parts to their connections, namely their tendencies to motion, or *conatus*, whose inherent actuality makes them appear a better candidate to ground physical substances. Moreover, having been defined as unassignable, it seemed natural to interpret the *conatus* as the physical counterpart of mathematical indivisibles, a topic that had captured the attention of Leibniz.

From the writings of this period, however, it is not clear how, and to what extent, the elimination of any material constituent can be finally implemented. Also the question of the coexistence of unassignable and assignable quantities is left unanswered. Probably unsatisfied by this incomplete picture, starting from the last years of his stay in Paris, here described as phase IV, Leibniz changes his approach: he dismisses the search for infinitesimal actual elements in physical bodies, and looks elsewhere to find the source of the actuality he needs.

More generally, he dismisses the search for any mathematical quantitative causes in physics, and turns to metaphysics, leaving to mathematics the (approximate) description of the phenomenal realm of nature, but not its causal explanation.

In the following subsection we analyze in more detail the notions of syncategorematicity and potentiality, and their relationship as conceived in the Scholastic logic.

### 2.2 Syncategorematicity and potentiality

The distinction between categorematic and syncategorematic infinite has a grammatical origin. C*ategoremata* were the parts of speech which possess a definite meaning of their own. *Syncategoremata* were such parts as depend on *categoremata* in order to acquire a definite meaning.[3]

Although in grammar *syncategoremata* and *categoremata* formed two mutually exclusive and collectively exhaustive classes, the notion of *syncategorema* was imported to logic, with a broadened range of application. For medieval logicians a *syncategorema* is a term that has a signification, but an indefinite one, and can be defined by the presence of other words.

In this context, a *categorema* can be employed in a "syncategorematic" way, and the distinction between categorematic and syncategorematic use of the same term becomes one of the standard ways of solving *sophismata*, namely puzzling sentences, that can be true or false depending on the interpretation of their parts.

The sentence "unlimited are limited" (*infinita sunt finita*) was a paradigmatic example of such a puzzle. The way of solving it by resorting to the opposition categorematic/syncategorematic made the distinction between categorematic and syncategorematic (use of the word) infinite a very common one.

Given the traditional Aristotelian distinction between actual and potential infinite, meant as philosophical notions, it is not surprising that the logico-grammatical level and the philosophical one became conflated, leading to talk about a categorematic and a syncategorematic infinite *tout court*. This is so, especially since the logicians themselves, introducing the grammatical distinction, referred to Aristotle's philosophical discussion of the infinite.

At least, this is the case of Peter of Spain, whose *Summulae Logicales* (Peter of Spain 2014), written in the second quarter of the thirteenth century, became the most influential introduction to logic of the late Middle Ages. Peter deals with the infinite in the last treatise of the book, devoted to distributions (*distributiones*). As we read in chapter 1:

> **[P1]** Distribution is the multiplication of a common term produced by a universal sign. When someone says 'every man,' for example, that term 'man' is distributed or confused for whatever is below it by this sign 'every,' and so in this case there is multiplication of a common term. [*Summulae Logicales* XII.1 Transl. Copenhaver][4]

After a short introduction, Peter discusses the most important distributive signs, and distinguishes them as *distributive of substance*: every, nothing, whole… and *distributive of accidents* (quality, quantity).

The infinite is the last among the signs distributive of quantity. Peter starts his discussion explaining the diverse meanings of the term "infinite", in philosophy. He quotes Aristotle's *Physics* almost verbatim:[5]

> **[P2]** Next comes the name 'infinite' which is said in five ways. In one way, the infinite is said to be what cannot be crossed-over – as a word is said to be infinite as to sight because it is invisible, since it is not naturally suited to be seen. In another way, the infinite is said to be what is incompletely crossed in that it has not yet been crossed, even though it is naturally suited to be crossed, as when someone is crossing a space and has not yet come to its limit.
>
> In a third way, the infinite is said to be in regard to succession, as an augmentable number is infinite when succeeded by unity or by another number. But in a fourth way it is in regard to division, like the continuous, for everything continuous is divisible without limit. Hence, it is

---

[3] Both the distinction and the term *syncategoremata* appear in the *Institutiones grammaticae*, written at the beginning of VI century by the Latin grammarian Priscian and well known in the Middle Ages. On syncategoremata in grammar and in logic see in particular (Kretzmann 1982).

[4] "Distributio est multiplicatio termini communis per signum universale facta. Ut cum dicitur omnis homo, iste terminus homo distribuitur sive confunditur pro quolibet suo inferiori per hoc signum omnis, et sic est ibi multiplicatio termini communis".

[5] *Physics* III 4, 204a2-7.

defined by Aristotle in the sixth book of the *Physics* as follows: "the continuum is divisible into what is always divisible."[6] And in a fifth way, the infinite is said in both those senses, by addition and by division, like time. For since time is continuous, it is divisible without limit, and in this way, it is infinite by division, and since one time comes after another time, in that way, by the addition of one time to another, it is infinite by addition. [*Summulae Logicales* XII.36 Transl. Copenhaver, with changes][7]

The first two ways of saying that something is infinite are in fact improper ways of speaking (as Aristotle himself explains), which Peter does not consider. He discusses the proper senses of the term, focusing on their crucial common feature: namely, their being defined in an iterative way:

> **[P3]** And as to these last three significations, the definition goes this way: the infinite is that whose quantity – for items that take quantity – can always get something more, so that if, after one part of a line, another part was taken on, and after that, a third, and the end of it could never be reached, then the line would be called 'infinite'. [*Summulae Logicales* XII.36 Transl. Copenhaver, with minor changes][8]

Here the reference is to Aristotle's definition of potential infinite:

> **[A]** the infinite is in virtue of another and another thing being taken, over and over again; and what is taken is finite, over and over again; but it is a different thing, over and over again [*Physics* III 6, 206a27–29. Translation ours].[9]

Then Peter comes to the dialectical parts of his discussion, and explains how infinite can be used both for a common and for a distributive term:

> **[P4]** But the usual claim is that the infinite is sometimes used for a common term, and then this proposition: "infinite are finite" is equipollent with this one "some infinite are finite".
> Sometimes, however, it is used for a distributive term, and then the former "infinite are finite" is equipollent as to distribution with this one: "the-more than-whatever are finite". And it is confirmed in this way: "the-more than-one are finite; the-more than-two are finite; the-more than-three are finite" and so on for others; "therefore, the-more than-whatever are finite".
> And then this is said to produce an 'inter-laddered' distribution (or else interrupted or discontinuous) because this word 'more' in the first proposition stands for two and beyond, in the second for three and beyond, and in this way, it keeps going higher step by step, or like a ladder. And therefore, the phrase 'more than whatever' produces an inter-laddered distribution because my term 'whatever' stands for some things and my term 'more' for others, as it goes higher

numerically—as has been said [*Summulae Logicales* XII.37 Transl. Copenhaver, with minor changes][10]

Finally, summarizing the solution of the sophisma, Peter maintains the equivalence between syncategorematic and distributive use of the term "infinite":

> **[P5]** A distinction made by others is that 'infinite' can be a common term, and in this way the original proposition is false; or else it can be a syncategorematic expression that in itself indicates distribution, as has been said, and in that way, they claim that the original proposition is true. [*Summulae Logicales* XII.38 Transl. Copenhaver, with minor changes][11]

In Peter's exposition the relation is clear, between the syncategorematic infinite, i.e., the infinite used for a distributive term, and the infinite in regard to addition, or to division, as described in **[P3]**: namely, Aristotle's potential infinite, as defined in **[A]**.

But while there is a general agreement concerning the identification of syncategorematicity with distributivity, scholars are divided on the possibility of equating syncategorematic and potential infinite.[12] This depends on the interpretation of Aristotle's definition of potential infinite: once the particular kind of potentiality Aristotle ascribes to the infinite is correctly understood, as in Peter's analysis, the identification with syncategorematicity is straightforward. On the contrary, no similarity can be detected between the two notions if one tries to ascribe to the infinite the "standard" notion of potentiality: namely, the one Aristotle usually deploys, and contrasts with actuality, when he speaks of objects, but explicitly prevents from using in the case of the infinite:

> We must not take potentiality here in the same way as that in which, if it is possible for this to be a statue, it actually will be a statue, and suppose that there is an infinite which will be in actuality [*Ph*. III 6, 206a18–20. Translation ours][13]

But then, what is the peculiar meaning Aristotle assigns to the notion of potentiality in the phrase "potential infinite"? The answer is in the definition itself:

> [A] the infinite is in virtue of (**A1**) another and another thing being taken, over and over again (*aei*); and (**A2**) what is taken is finite, over and over again (*aei*); but (**A3**) it is a different thing, over and over again (*aei*). [*Physics* III 6, 206a27–29. Translation ours][14]

As is clear, and as was clear to Peter of Spain, **[A]** is not a conventional Aristotelian definition, but an operational characterization. It does not describe an object but an action – 'taking' something – and the iteration of the action, signaled by the adverb *aei*. For each thing we take (*first step*), there is another thing to take (*second step*), and another to take beyond that (*third step*), and another beyond, and so on, over and over again (*aei, next steps*).

Once the potential infinite is understood in this way, it is immediate to interpret it as "a syncategorematic expression that in itself indicates distribution" to use Peter's words.[15] In fact, as he explains in **[P4]**, if the term "infinite" is used in this sense, the sentence "infinite are finite" can be rephrased as: "the-more than-one are finite (*first step*); the-more than-two are finite (*second step*); the-more than-three are finite (*third step*) and so on for others (et sic de aliis, *next steps*)".[16]

### 3. LEIBNIZ'S TRAJECTORY IN FIVE PHASES

*3.1 An overview*

Leibniz's mature position concerning the ultimate composition of bodies, and the related issue of their division, is rather clear (see Phase V, below). It is part of an articulated philosophical construction, embracing mathematics, physics, and metaphysics, that can be summarized as follows.

Since a complete separation exists between the ideal realm of mathematics, the phenomenal realm of physics, and the true realm of metaphysics, a completely different solution is offered to the problem in the three contexts.

Metaphysical real entities are indivisible unextended unitary wholes, which do not admit division of any kind. However, these substances, as any other metaphysical notion, are unknowable and inaccessible per se. As unknowable, in particular, and extraneous to the category of quantity, they can be described only through analogical speech; as inaccessible, they can be approached only through their phenomenal "expressions", namely physical bodies.[17]

Physical bodies are extended and composed of extended parts. These parts are actually infinite, but only in this sense: there is no part that is not actually subdivided into others:

> It is perfectly correct to say that there is an infinity of things, i.e., that there are always more than one can specify. But it is easy to demonstrate that there is no infinite number, nor any infinite line or other infinite quantity, if these are taken to be genuine wholes. The Scholastics were taking that view, or should have been doing so, when they allowed a 'syncategorematic' infinite, as they called it, but not a 'categorematic' one. [*Nouveaux essais sur l'entendement humain*, II. XVII. A6.6 157. Transl. Remnant-Bennett].[18]

---

[14]  See note 9 for the Greek text. The condition **A3** has the only aim of discriminating between procedures that are infinite in a proper sense and procedures that are infinite because they are periodic. In the latter case, one comes back to something already taken after a certain number of steps

[15]  See also William of Sherwood's *Syncategoremata*, where the infinite in its syncategorematic significance is explicitely related to induction (O'Donnell 54-55, translated in Kretzmann 42-43).

[16]  The sophisma itself "infinite are finite" is nothing but a deliberately equivocal way of restating Aristotle's condition **A2**.

[17]  The notion of "expression" is repeatedly discussed throughout the correspondence with Arnauld from September 1687 [A2.2. 230-8].

[18]  "A proprement parler, il est vrai qu'il y a une infinité de choses, c'est à dire qu'il y en a tousjours plus qu'on n'en peut assigner. Mais il n'y a point de nombre infini ny de ligne ou autre quantité infinie, si on les prend pour des veritables Touts, comme il est aisé de demonstrer. Les écoles ont voulu ou dû dire cela, en admettant un infini syncategorematique, comme elles parlent, et non pas l'infini categorematique".

Moreover, the parts of a physical body do not form a continuous magnitude, but only a contiguous one. Here one must note that neither the notion of actual infinity employed in physics – also known in literature as "syncategorematic actual" – nor that of contiguity are mathematical notions.[19]

In fact, in mathematics Leibniz maintains his traditional genuinely syncategorematic characterisation of continuity as potentially infinite divisibility, which means that the parts into which a geometrical continuum can be divided are only potential:

> As for the first point, it follows from the very fact that a mathematical body cannot be analyzed into primary constituents that it is also not real but something mental and designates nothing but the possibility of parts, not something actual. A mathematical line, namely, is in this respect like arithmetical unity; in both cases the parts are only possible and completely indefinite. A line is no more an aggregate of the lines into which it can be cut than unity is the aggregate of the fractions into which it can be split up. […]
> But in real things, that is, bodies, the parts are not indefinite – as they are in space, which is a mental thing – but actually specified in a fixed way according to the divisions and subdivisions which nature actually introduces through the varieties of motion. And granted that these divisions proceed to infinity, they are nevertheless the result of fixed primary constituents or real unities, though infinite in number. Accurately speaking, however, matter is not composed of these constitutive unities but results from them. [*Leibniz to De Volder* June 30 1704; G II, 268. Transl. Loemker][20]

Despite the differences, both in physics and in mathematics no infinitesimal unassignable quantities can be reached while proceeding with the division: the quantity obtained at any step, in fact, is in a given and finite ratio with the one obtained at the previous step, and with the one obtained at the next one.

But while to postulate the existence of infinitesimal quantities in physics could raise some difficulties,[21] the introduction of infinitesimal quantities in mathematics does not carry such risks, and is highly useful. Although they do not follow Euclid V Definition 4 and do not admit a traditional geometrical representation, Leibniz provides rules that make his unassignable infinitesimals mathematically sound, non-contradictory objects.[22]

---

[19]  On the expression "syncategorematic actual" for denoting Leibniz's physical infinite, which is currently employed in Leibnizian studies, see Arthur (2015, 145). It has been defended as an innovative mathematical notion for instance in Arthur 2019, criticized on that point in Bussotti 2019. On the mathematical inconsistency of Leibniz's physical notion of infinite see also Ugaglia 2022.

[20]  "Quod primum attinet, eo ipso quod corpus mathematicum non potest resolvi in prima constitutiva, id utique non esse reale colligatur, sed mentale quiddam nec aliud designans quam possibilitatem partium, non aliquid actuale. Nempe linea mathematica se habet ut unitas arithmetica, et utrobique partes non sunt nisi possibiles et prorsus indefinitae; et non magis linea est aggregatum linearum in quas secari potest, quam unitas est aggregatum fractionum in quas potest discerpi […] At in realibus, nempe corporibus, partes non sunt indefinitae (ut in spatio, re mentali), sed actu assignatae certo modo, prout natura divisiones et subdivisiones actu secundum motuum varietates instituit, et licet eae divisiones procedant in infinitum, non ideo tamen minus omnia resultant ex certis primis constitutivis seu unitatibus realibus, sed numero infinitis. Accurate autem loquendo materia non componitur ex unitatibus constitutivis, sed ex iis resultat".

[21]  Even if the existence of infinitesimal quantities in physics is not an absolute impossibility – namely, it does not involve any contradiction – it is contrary to the harmony of things. In fact, a world containing "some portion of matter which has an unassignable ratio to another portion" would be less perfect than a world where all portions are comparable. On this point it is interesting to read the 1698 correspondence between Leibniz and Bernoulli, vigorous proponent of the existence of infinitesimal quantities in physics. In particular, in his letter of July 5 1698 (GM III 503-504), Bernoulli asks why only a size (*gradus*) of infinity is allowed, namely the one proportionate to our intellect, and not the infinite others (he means here infinitesimals and higher degree of infinitesimals dx, ddx...). In his answer of July 12, Leibniz replies that: "from the actual division it follows that in every part of matter, there is a world – so to speak – made of innumerable creatures; but we ask if there is some portion of matter which has an unassignable ratio to another portion, or if there is a line bounded (terminata) on both ends, having an infinite or an infinitely small ratio, to another line. We assume that in our Calculus, fruitfully, but it does not follow that this holds in nature. But this is the subject of another research" (GM III 516). Accidental impossibility of Leibniz's infinitesimals has been argued in Esquisabel-Quintana 2021.

In other words, while geometrical curves are continuous, namely indefinitely divisible into potential parts, the problems they represent can be solved by using infinitesimal quantities. Moreover, the language of infinitesimals turns out to be useful in speaking, analogically, of some metaphysical objects.

The situation is similar to that with imaginary numbers, which at Leibniz's time had neither a geometrical representation, nor a physical counterpart, but were well defined, non-contradictory and useful mathematical notions. And in fact, as for imaginary numbers, in his mature philosophy Leibniz refers to infinitesimals using the language of imagination – to feign, fictions – and contrasts them with more traditional mathematical notions, referred to as ideal concepts, which are the image, although idealized, of physical situations.[23]

### 3.2 Five phases

As mentioned above, Leibniz reached this final position through a complex path, which we will trace in the present section. We have focused our attention on five phases, which are especially significant for reconstructing the entire evolution. Except for Phase I, where Leibniz's physical and mathematical analysis do coincide, we have treated them separately, as far as possible.

### Phase I (1668-9)

While in his mature philosophy Leibniz will analyze physical bodies in terms of individual metaphysical substances, he was initially confident that one should banish from physics any metaphysical explanations, and replace them by purely mechanical ones. In particular, Leibniz attempted to reduce the Aristotelian triad *matter-form-change*,[24] which he considered too abstract, to the "mechanical" triad *extension-shape-motion*. Wrote Leibniz:

> The one question is whether Aristotle's abstract theories of matter, form, and change can be explained by magnitude, shape, and motion. This the Scholastics deny and the philosophical reformers affirm. The latter opinion seems to me to be not only the more true but also the more consistent with Aristotle. [*Leibniz to Thomasius*, April 1669. A2.1². 25-26. Transl. Loemker with minor changes][25]

A physical body is the compound of a passive element, or matter, which Leibniz identifies with extension,[26] and an active one, or form, which he identifies with shape. Both extension and shape are dealt with in a rather simplistic way, in accordance with the fact that in this phase Leibniz's knowledge of mathematics, which essentially means Greek geometry, is rather basic, and for the most part acquired by reading philosophical works. In particular, extension is continuous, and continuity means for Leibniz potentially infinite divisibility, that is syncategorematicity, as it was for Aristotle.[27] The terminology itself of unassignability and infinitesimals is absent in this phase: "infinitely small" just means "potentially infinitely divisible".[28]

---

22  The rules are discussed in the *Nova Methodus pro maximis et minimis*, published in 1684 in the *Acta Eruditorum* [GM V pp. 220-226].

23  On Leibniz's ideal notions see Sherry-Katz 2014 and Katz-Sherry 2013.

24  The noun *change* is our rendering of Aristotle's μεταβολή, which includes motion (κίνησις, namely local change) as well as qualitative and quantitative changes.

25  "Hoc unum in quaestione est, an quae Aristoteles de materia, forma, et mutatione abstracte disputaverit, ea explicanda sint per magnitudinem, figuram, et motum. Id Scholastici negant, Reformatores affirmant. Reformatorum sententia mihi non solum verior, sed et Aristoteli magis consentanea videtur". Here and throghout the paper we have rendered "figura" with "shape", instead of Loemker's "figure".

26  In physical bodies extension is coupled with impenetrability, which sometimes Leibniz prefers to call antitypy (ἀντιτυπία): "To be extended, however, is nothing else but to be in space, and antitypy is the impossibility of being in the same space with another thing" [*Leibniz to Thomasius*, April 1669. A2.1². 36. Transl. Loemker].

> Matter has quantity too, though it is indefinite, or interminate as the Averroists call it. For being continuous, it is not cut into parts and therefore does not actually have boundaries. [*Leibniz to Thomasius*, April 1669. A2.1². 26-27. Transl. Loemker][29]
> Space is a primary extended being or a mathematical body, which contains nothing but three dimensions and is the universal locus of all things. Matter is a secondary extended being, or that which has, in addition to extension or mathematical body, also a physical body, that is, resistance, antitypy, solidity, the property of filling space, impenetrability[30] [*Leibniz to Thomasius*, April 1669. A2.1². 34. Transl. Loemker with minor changes][31]

Motion intervenes in the generation of physical bodies. In particular, it is responsible for their shape:

> Form (in the technical sense) as educed from the potentiality of matter, is nothing more than this particular shape of the whole as generated from this particular motion of matter, and from this particular situation of parts: e.g., the form of a square $\boxdot c$ is educed from the union of two triangles, $\Delta a$ and $\nabla b$, and from the potentiality of their matter to be in reciprocal contact by motion. [*Leibniz to Thomasius*, October 1668. A2.1².17][32]

The world is a plenum, made of continuous bodies, actually divided from each other:

> This makes it clear that, if mass were created discontinuous or separated by emptiness in the beginning, there would at once be certain concrete forms of matter. But, if it is continuous in the beginning, forms must necessarily arise through motion. For division comes from motion, the bounding of parts comes from division, their shapes come from this bounding, and forms from shapes; therefore, forms come from motion […] From this it is clear that every arrangement into

---

[27] Aristotle characterizes the continuum as "infinitely divisible" (*Physics* III 1, 200b18; I 2, 185b10–11; VI 6, 237a33; VI 8, 239a22) or "divisible into what is always (*aei*) divisible" (*Physics* VI 2, 232b24–25; see IV 12, 220a30; VI 1, 231b15–16; VI 6, 237b21; VIII 5, 257a33–34; *On the heavens* I 1, 268a6–7). Infinite divisibility is an essential property of the *continuum*, which directly follows from Aristotle's definition of continuity, as a relational notion: two parts of a whole are said to be continuous with one another when they are in contact and their ends therefore become one (συνεχῆ μὲν ὧντὰ ἔσχατα ἔν *Physics* VI 1, 231a21-22 "we defined continuous things whose ends are one"); see also *Physics* V 4, 228a29-30; V 3, 227a11-12; *Cat.* 6, 4b20-5a14. We discussed the continuum as always (*aei*) divisible in Subsection 2.2 above.

[28] Cf. the analysis of the infinitely small in *Physics* III.6, where Aristotle argues against the existence of indivisible lines (ἀτόμους γραμμάς 206a17-18).

[29] "Quantitatem quoque habet materia, sed interminatam, ut vocant Averroistae, seu indefinitam, dum enim continua est, in partes secta non est, ergo nec termini in ea actu dantur".

[30] Resistance, antitypy, solidity, the property of filling space and impenetrability are synonyms, and refer to the impossibility, for a physical body, of being in the same space with another thing, as Leibniz explains a few lines later.

[31] "Spatium est Ens primo-extensum, seu corpus mathematicum, quod scilicet nihil aliud continet quam tres dimensiones, estque locus ille universalis omnium rerum. Materia est ens secundo-extensum, seu quod praeter extensionem vel corpus mathematicum habet et corpus physicum, id est, resistentiam, ἀντιτυπίαν, crassitiem, spatii-repletivitatem, impenetrabilitatem".

[32] "Formam enim educi ex potentia materiae, nihil aliud est, quam ex hoc materiae motu, ex hoc partium situ, hanc totius figuram oriri: v.g. ex unione duorum triangulorum, $\Delta a$ et $\nabla b$, et ex potentia materiae eorum ad contactum mutuum per motum, educitur forma quadrati $\boxdot c$." In this case, the resulting square does not keep trace of the way in which it has been created. The same example of the square arising from two triangles appears in a 1676 note – see footnote 71 – but here Leibniz resorts to metaphysics in order to obtain a form which keeps trace of its constituents.

a form is motion, and the vexatus problem of the origin of forms is answered.[33] [*Leibniz to Thomasius*, April 1669 A2.1².26-27. Transl. Loemker][34]

Phase I in summary:
Physics: physical objects are continuous; continuous = in(de)finitely divisible; parts are only potential and finite. Both the approach and the terminology employed are syncategorematic.
Mathematics: mathematical objects are continuous; continuous = in(de)finitely divisible; parts are only potential and finite. Both the approach and the terminology employed are syncategorematic.

## Phase II (1670-1)

Continuity as conceived by Aristotle, and as employed by Leibniz in phase I, means potentiality of the parts, and potentiality means indefiniteness. However, indefiniteness, in Leibniz's opinion, has more to do with Descartes' thoughts that with reality; see item (2) below. Phase II represents Leibniz's first attempt to overcome this difficulty without abandoning his initial purpose of leaving aside any metaphysics. However, Leibniz's poor knowledge of mathematics will prove to be a hindrance to effective implementation of the project.[35]

*Phase II in physics*
The project was to give some sort of individuality to the parts of a physical body, still treating them as continuous. For doing so, Leibniz broadens the notion of continuity, and affirms – without any concern for the mathematical consistency of the statement – that a continuous body is not merely infinitely divisible, but it is infinitely divided into parts, and these parts are actual:

(1) *There are actually parts in the continuum*
(2) *these are actually infinite*, for Descartes's indefinite is not in the thing but in the thinker [*Theoria motus abstracti* 1671, A6.2. 264. Transl. Arthur][36]

Being actual, the parts of a physical body are organized in a precise way. Therefore, in this phase form does not depend only on the external shape of the body, as in phase I, but also on the disposition of its internal parts:

---

[33] As Leibniz himself wrote some years later: "Tempus erat quo credebam, omnia Motuum Phaenomena principiis pure Geometricis explicari posse, nulli Metaphysicis propositionibus assumptis" [There was a time when I believed that all the phenomena of motion could be explained on purely geometrical principles, assuming no metaphysical propositions. *Principia mechanica ex metaphysicis dependere* 1678/81 A6.4.1976, Transl. Ariew-Garber].

[34] "Ex his patet siquidem ab initio massa discontinua seu vacuitatibus interrupta creata sit, formas aliquas statim materiae concreatas esse; sin vero ab initio continua est, necesse est, ut formae oriantur per motum (nam de annihilatione certarum partium ad vacuitates in materia procurandas, quia supra naturam est, non loquor) quia a motu divisio, a divisione termini partium, a terminis partium figurae earum, a figura formae, ergo a motu formae. Ex quo patet, omnem dispositionem ad formam esse motum, patet quoque solutio vexatae de origine formarum controversiae."

[35] On Leibniz's knowledge of mathematics before his stay in Paris see for instance Hofmann 1974, pp. 1-11.

[36] "Dantur actu partes in continuo […] eaeque infinitae actu, indefinitum enim Cartesii non in re est, sed in cogitante". We find the same assumption in Leibniz's *Hypothesis Phisica Nova*, dated shortly before and also known as *Theory of concrete motion*: "Quae res ob divisibilitatem cujuslibet continui in partes quantumvis parvas in infinitum non est difficilis explicatu" [*Hypothesis Phisica Nova* 1671. A6.2 224]. Here physical bodies are thought as made of parts, or *bullae* (bubbles), in a more or less solid state, depending on the circumstances: "Principio autem ex fluidi aestuatione et fusione per lucem seu calorem ortae sunt bullae innumerabiles ac magnitudine crassitieque variantes. Nam quoties subtilia perrumpere per densa conantur, et est quod obsistat, formantur densa in cavas quasdam bullas, motumque partium internum, ac proinde consistentiam seu cohaesionem (per nostram *de motu Theoriam*, theor. 17) nanciscuntur." [*Hypothesis Phisica Nova* 1671. A6.2 226].

For Aristotle too, substantial form (with the exception of mind) is not an absolute entity, but only a *logos*, a ratio, a proportion, a number, a structure of internal parts. [*Hypothesis Phisica Nova* 1671. A6.2. 247].[37]

But how can these parts stay together, forming a continuous body? To solve the problem Leibniz resorts to motion: being in motion, parts experience an endeavor (*conatus*)[38] to penetrate each other, and this causes their extremities to become one. But the coincidence of the extremities of the parts is the condition Aristotle requires in order to have a continuous, unitary whole (see note 27):

> I should think that the *conatus* of the parts toward each other, or the motion through which they press upon each other, would itself suffice to explain the cohesion of bodies. For bodies which press upon each other are in a *conatus* to penetrate each other. The *conatus* is the beginning (*initio*); the penetration is the union. But when bodies begin to unite, their beginnings or extremities (*termini*), are one. Bodies whose extremities are one, or τὰ ἔσχατα ἕν, are according to Aristotle's definition not only contiguous but continuous, and truly one body, movable in one motion [*Leibniz to Hobbes*, July 1670. A2.1². 92. Transl. Loemker with minor changes][39]
> If one body strives against another, this is the beginning of mutual penetration or union or that the boundaries of the two are one, as Aristotle defines a continuum [*Leibniz to Arnauld*, November 1671. A2.1². 279. Transl. Loemker][40]

In his more applied writings Leibniz describes this situation in concrete terms. The behavior of the parts of physical bodies, which he calls "bubbles", is explained in analogy with what happens in a turbulent fluid: as the bubbles which are continuously produced and destroyed in a fluid, also the parts of a body are elastic and constantly move, breaking and recombining with each other.

*Phase II in mathematics*
As in phase I, Leibniz does not draw a clear distinction between physics and its mathematical description, but while in phase I continuity as potentially infinite divisibility had a clear meaning both in physics and in mathematics – a body is made of continuous parts in the same way as a segment is made of segments – the composition of bodies as conceived by Leibniz in phase II is more difficult to translate in mathematical terms.

The main problem is the mathematical *status* of the actually infinite parts, and their relationship with unassignable endeavors.

Later in the *Theory of abstract motion*, we learn that Leibniz rejects *minima* (i.e., points as constituent parts of the continuum in the modern sense) and *maxima*, where a maximum is an unbounded infinite line (the kind we have in modern Euclidean geometry). The reason for both is a violation of the part-whole principle that would ensue via Galileo's paradox:

---

[37] "Certe formas substantiales (demta mente) etiam Aristoteli non esse ens absolutum, sed tantum λόγον, rationem, proportionem, ἀριθμόν, structuram partium intimam".

[38] Leibniz borrows the notion of *conatus* from Hobbes' endeavour, which has a very similar definition, and a very similar role in physics. However, while Leibniz's *conatus* are unextended, Hobbes' endeavours are finite quantities, although negligible, like his points: "I define endeavour to be Motion made in less Space and Time then can be given; that is, less than can be determined or assigned by Exposition or Number; that is, Motion made through the length of a Point, and in an Instant or Point of Time" [*De Corpore* 3.15.2]. For Hobbes a point is not "that which has no quantity" but "that whose quantity is not at all considered" [*ibid.*]. On the role of Hobbes as the mediator between the young Leibniz and mathematics, and particularly Cavalieri's method, see for instance Hofmann 1974.

[39] "Ego crediderim ad cohaesionem corporum efficiendam sufficere partium conatum ad se invicem, seu motum quo una aliam premit. Quia quae se premunt sunt in conatu penetrationis. Conatus est initium, penetratio unio. Sunt ergo in initio unionis. Quae autem sunt in initio unionis, eorum initia vel termini sunt unum. Quorum Termini sunt unum seu τἀἔσχαταἔν, ea etiam Aristotele definitore non jam contigua tantum, sed continua sunt, et vere unum corpus, uno motu mobile."

[40] "Si corpus conetur in corpus, esse ambo in initio penetrationis seu unionis, seu extrema eorum unum esse ut continuum definit Aristoteles [ὧν] τὰ ἔσχατα ἕν."

(3) *There is no minimum in space or body*, that is, nothing which has no magnitude or part. For such a thing has no situation, since whatever is situated somewhere can be touched by several things simultaneously that are not touching each other, and would thus have several faces; nor can a minimum be supposed without it following that the whole has as many minima as the part, which implies a contradiction [*Theoria motus abstracti* 1671, A6.2. 264. Transl. Arthur].[41]

On the other hand, he acknowledges the existence of indivisibles – namely, the mathematical counterpart of endeavors – and argues as follows:

(4) *There are indivisibles or unextended things*, otherwise neither the beginning nor the end of a motion or of a body is intelligible. This is the demonstration: any space, body, motion, and time has a beginning and an end. Let that whose beginning is sought be represented by the line *ab*, whose midpoint is *c*, and let the midpoint of *ac* be *d*, that of *ad* be *e*, and so on. Let the beginning be sought to the left, on *a*'s side. I say that *ac* is not the beginning, since *dc* can be taken away from it without destroying the beginning; nor is *ad*, since *ed* can be taken away, and so on. Therefore, nothing is a beginning from which something on the right can be taken away. But that from which nothing having extension can be taken away is unextended. Therefore, the beginning of a body, space, motion, or time (namely, a point, an endeavour, or an instant) is either nothing, which is absurd, or is unextended, which was to be demonstrated [*Theoria motus abstracti* A6.2. 264. Transl. Arthur with minor changes][42]

Endeavors are unextended amounts of motion, and indivisibles, or "fat points" (our term) are unextended portions of bounded lines.

Moreover, unlike Euclidean or Aristotelian points, which are "unextended beginnings" of lines but do not enter in their composition, Leibnizian "fat points" do make lines (and "fat lines" do make surfaces...).

As a consequence, even if they are smaller than any assignable quantity, "fat points" can be compared with each other. In Leibniz's words, they are different in "magnitude": bigger in a bigger curve, smaller in a smaller one:

The cone is made by as many circles parallel to the basis as points of the height [*Nugae pueriles* 1670. A7.4. 57-58][43]

Whence the unassignable arc of a bigger circle is greater than that of a smaller one [*Theoria motus abstracti* A6.2. 267. Transl. Arthur][44]

An *arc* smaller than any that can be given is still greater than its chord, although this is also smaller than can be expressed, i.e., consists in a point. But that being so, you will say, an *infinitangular polygon* will not be equal to a circle: I reply, it is not of an equal magnitude, even if it be of an equal extension: for the difference is smaller than can be expressed by any number [*Theoria motus abstracti* A6.2. 267. Transl. Arthur][45]

A point can be bigger than another, but in a smaller ratio (*in ratione minore*), than the one that can be exhibited, or infinite with respect to any sensible quantity [*Leibniz to Arnauld*, November 1671. A2.1². 278][46]

The main problem with this picture is the necessity of distinguishing between "having extension" and "having magnitude", a problem that Leibniz thinks he has solved by employing the locution "having magnitude" as synonymous with "having parts", where parts are not necessarily "external to each other". In fact, Leibniz says, parts can be indistant, and this is the basis of Cavalieri's method too:

(5) *point is not that which has no part*, nor that whose part is not considered; but that which *has no extension*, i.e., whose parts are indistant, whose magnitude is inconsiderable, unassignable, is smaller than can be expressed by a ratio to another sensible magnitude unless the ratio is infinite, smaller than any ratio that can be given. But this is the basis of the *Cavalierian Method*, whereby its truth is evidently demonstrated, inasmuch as one considers certain rudiments, so to speak, or beginnings, of lines and figures smaller than any that can be given. [*Theoria motus abstracti* A6.2. 265. Transl. Arthur][47]

The nature of points is undeniably admirable: indeed, although a point is not divisible into parts external to each other, however it is divisible into parts which before were not external to each other, or into parts that before were penetrating each other. [*Leibniz to Oldenbourgh*, September 1670. A2.1². 103][48]

Summarizing, Leibniz tries to reconcile the actuality of the infinite parts into which physical bodies are divided – and not only divisible – and the potentiality of matter. As material parts, in fact, their magnitude can decrease ad infinitum (no *minima*).[49]

But the conflation between actuality and potentiality is not the only flaw in Leibniz's construction. Furthermore, the physical notion of *conatus*, that is of an endeavor of motion, creates problems when one tries to translate it into mathematics.

Two physical endeavors can have the same extension (in fact, they have no extension) but a different "intension": for instance, when one begins a faster motion than the other. But what does it mean that two unextended "fat points" do have a different "magnitude"? Leibniz's explanation in term of indistant inner parts, which can penetrate each other, does not have a clear meaning in mathematics. However, this analysis is never elaborated by Leibniz, who merely claims that with the identification of his "fat points" with indivisibles he has demonstrated Cavalieri's method, too.[50]

---

[45] "Arcus minor quam qui dari potest, utique chorda sua major est, quanquam haec quoque sit minor quam quae exponi potest, seu consistat in puncto. At ita, inquies, Polygonum infinitangulum non erit circulo aequale: respondeo, non esse aequalis magnitudinis, etsi sit aequalis extensionis: differentia enim minor est quam ut ullo numero exprimi possit"

[46] "Esse punctum puncto majus, sed in ratione minore, quam quae exponi potest, seu ad sensibilem quamcumque infinita".

[47] "Punctum non est, cujus pars nulla est, nec cujus pars non consideratur; sed cujus extensio nulla est, seu cujus partes sunt indistantes, cujus magnitudo est inconsiderabilis, inassignabilis, minor quam quae ratione, nisi infinita, ad aliam sensibilem exponi possit, minor quam quae dari potest: atque hoc est fundamentum Methodi Cavalierianae, quo ejus veritas evidenter demonstratur, ut cogitentur quaedam ut sic dicam rudimenta seu initia linearum figurarumque qualibet dabili minora"

[48] "Admirabilis profecto est natura punctorum; quanquam enim punctum non sit divisibile in partes positas extra partes, est tamen divisibile in partes antea non positas extra partes, seu in partes antea se penetrantes"

[49] A different interpretation has been propounded in Arthur (2009, pp. 17-18), which identifies actually infinite parts with points: in this case the whole physical construction could be traced back to Cavalieri's theory, and Leibniz's construction would appear more consistent. In principle this would be an appealing interpretation, the only problem being that there is no basis for it in Leibniz's text. In particular, he never identifies the physical notions of part and point, and repeatedly refers to points as extremities of parts.

Phase II in summary:
Physics: physical objects are actually infinitely divided into continuous material parts (no *minima*) whose extremities coincide. Extremities are held together by unextended endeavors of motion.
Mathematics: material extension is represented by a continuum, as in phase I. Endeavors of motion are represented by unextended "fat points", unassignable and indivisible, which have parts but parts that are indistant.

## Phase III (1672-4)

Phase III, which coincides with the first years of Leibniz's stay in Paris, marks the beginning of his mathematical career. It is in Paris that Leibniz undertakes serious training in mathematics, and enters in relation with Huygens, and then with other leading mathematicians, who will become his correspondents.

In phase II, Leibniz's knowledge of the contemporary methods for solving the problems related to infinity was rather poor. He realized the importance of Cavalieri's indivisibles, and more generally of unassignable infinitely small quantities, but he lacked the technical tools necessary in order to transform a philosophical intuition into a mathematical investigation.

The idea of a point with indistant parts, for instance, while appealing and maybe effective in a dialectical argument, was very hard to disambiguate, and to translate into a sound mathematical notion. Analogously, the coexistence of actuality and potentiality implied in Leibniz's analysis of physical bodies, would lead to serious problems, if mathematically analyzed.

But in phase II Leibniz did not seem aware of the existence of such difficulties, which become clear only in subsequent years, as his mathematical knowledge grew. For this reason, even if in terms of effective results phase III is merely a transitional phase, it is a crucial step in our analysis. Even though Leibniz did not publish anything, the large amount of handwritten notes he left, often inconclusive and vitiated by inconsistencies, show nevertheless the awareness of the impasse which phase II was leading to, and the attempt to avoid it.

### Phase III in physics

The first novelty, which Leibniz introduces in physics in order to keep actuality and potentiality separate, consists in reversing the arrow of the analysis: if, by dividing a continuum, one can only reach potentially divisible parts, the only way to obtain a matter made of genuine actual ultimate parts is to start from these parts.

Both in phases I and II, the basic stuff – from which physical bodies was derived – was a continuous entity: once continuity was broken, in virtue of motions, parts were created, which entered into the composition of physical bodies. In phase III, by contrast, the ultimate element is a plurality of basic material entities, discrete from the outset, spherical, which aggregate in order to give rise to physical bodies. Leibniz calls them "atoms" or *bullae*, or *terrellae*:

> Sensible bodies are aggregates of countless (innumerabilium) terrellas [*Propositiones quaedam physicae*, 1672. A6.3.32][51]

---

[50] On the ambiguous relation between indivisibles and motion see Cavalieri's reaction to Guldin's critique (*Exercitationes geometricae sex*, Exercitatio III, cap. VII). Cf. Cavalieri 1966 tr. Radice 808-810.

[51] "Corpora sensibilia sunt aggregata innumerabilium terrellarum". A *terrella* is a little spherical loadstone, as described in Gilbert, *De Magnete* 1600.

By introducing discrete elements, which have a proper individuality and self-sufficiency, Leibniz solves any ambiguity concerning their *status*: the ultimate parts are actual and indivisible.

The problem is that these atoms are meant to be part of bodies which are continuous: how is this possible? Leibniz solves the problem in one direction, showing how a continuous stuff can be obtained starting from discrete elements. As in phase II, in fact, he resorts to motion for welding together the material atoms which form a body: due to their endeavors to motion, the discontinuities between different atoms disappear:

> If a body is forced (*conetur*) into the place of another, these two bodies are continuous [*De consistentia corporum* 1672/3. A6.3.95][52]

Moreover, due to their smallness, even the empty spaces which exist between the spherical atoms become negligible, and do not interrupt continuity:

> It is necessary that matter is uninterrupted, and it completely fills the space, or at least it leaves some unnoticeable (*insensibiles*) empty islands, without any breaking of continuity. [*Demonstratio substantiarum incorporearum*, 1672. A6.3. 78][53]

But once physical bodies are continuous, they must be also infinitely divisible:

> This sounds paradoxical - I must admit - but is no less necessary than the actual division of the continuum, in infinite parts (which follow from this exact proof, consistently with Cartesian teachings) and other marvels of the internal nature [*Propositiones quaedam physicae*, 1672. A6.3.57][54]

And then the question arises, how can individuality and infinite divisibility coexist?

Reading Leibniz's notes, it seems that the only possible answer is that they cannot, and that the only way to obtain a consistent description of physical bodies is to eliminate any (atom of) matter:

> From this we understand that matter and motion, or *conatus*, are the same […] And getting back to our issue, it is clear from these facts that matter does not consist of extension, but of motion, and a lot of matter is nothing more than a lot of motion. [*Propositiones quaedam physicae*, 1672. A6.3 56][55]
> Hence it follows that there is no matter in body distinct from motion, since it would necessarily contain indivisibles [*De minimo et maximo* 1672/3. A6.3. 100][56].

In other words, motion can contain indivisibles, matter cannot. And here two difficulties arise: is a motion made of unassignable indivisibles continuous? And if it is, in what sense? Moreover, if bodies are nothing but motion, what happened to the atoms of matter? Leibniz's writings of this period do not contain an answer, and in fact, even if the starting point was different, the physical theory of continuity Leibniz builds in phase III ends up with the same inconsistencies as the one he built in phase II.

---

[52] "Si corpus unum conetur in alterius locum, ea duo corpora continua sunt".

[53] Leibniz will come back on the status of *vacua intermixta*, which he calls also *vacua metaphysica* and their incomparable magnitude some years later, in phase IV: "[…] because spheres, however many they may be, do not fill a whole. However, all emptinesses collected into one have no greater ratio to any assignable space than the angle of contact has to a straight line" [*De summa rerum* 1675/6. A6.3. 526. Transl. Parkinson]. "A metaphysical vacuum is an empty place, no matter how small, provided that it is genuine and real. A physical plenum is consistent with a metaphysical vacuum which is unassignable" [*De summa rerum* 1675/6. A6.3. 473. Transl. Parkinson]. Note that in this period and in this context the term "metaphysical" still retains a negative connotation, keeping track of the "empty words" imputed to Scholasticism.

[54] "Res est fateor paradoxa sed non minus necessaria, quam divisio continui actualis, in partes infinitas (quae ex hac ipsa demonstratione confessione Cartesii sequitur) aliaque naturae interioris admiranda"

[55] "Unde intelligi potest, materiam et motum vel conatum idem esse […]. Patet ergo ex his, ut viam redeam, materiam non consistere in extensione, sed in motu, et nihil aliud esse mulatm materiam, quam multum motum".

[56] "Hinc sequitur nullam esse in corpore materiam a motu distinctam, ea enim necessario indivisibilia contineret".

*Phase III in mathematics*
In mathematics the situation is less problematic, although some doubt remains concerning the composition of the continuum. Apparently, Leibniz's mathematical continua can be read as entirely composed of "fat points", since now these points are no longer seen as unextended, but as extended infinitely small lines. As their ancestors, these new points are unassignable, i.e., infinitely smaller than any assignable line:

> Continuous space is made of parts smaller than anything determinable by us [*Demonstratio substantiarum incorporearum*, 1672. A6.3.81]
> In a continuum there are some infinitely small things, or infinitely smaller than any perceptible one [*De minimo et maximo* 1672/3. A6.3.98].
> A point can be infinitely smaller than another [*De minimo et maximo* 1672/3. A6.3.99][57]

In this period "infinitely small", "unassignable" and "point" are employed as synonyms:

> In short: each of these scalene triangles is a point, because their base and height are infinitely small lines. For this reason, all together, arranged on a single straight line obtained from the extended evolute (in fact, the triangles are as many as the points of the chord) do not make but a line [*Varia de cycloeide* 1673 A7.4 p. 74][58]

Here indivisibles are defined as synonymous with infinitely small:

> Indivisibles are to be defined as infinitely small, or as having infinite ratio, or difference, with a sensible quantity [A7.4 p. 265][59]

During these years the term "infinitesimal" makes its appearance, but it is always employed in the phrase "infinitesimal (part) of", and denotes a particular case of infinitely small line: the infinitesimal part of a quantity is one of the infinitely many equal parts into which the quantity has been divided:[60]

> We suppose that the abscissa AE is divided into equal infinitely many equal parts as EF, FG […]
> And the straight line EF, or FG will be infinitely small, and each one will be an infinitesimal of the abscissa AE [*Methodus nova investigandi tangentes linearum curvarum...* 1673 A7.4 p. 657][61]
> And this is the proof: suppose that MNE is the characteristic triangle of the figure, that is the one whose height MN is an infinitesimal of the height of the figure [*Trigonometria inassignabilium* 1673 A7.4 p. 497][62]

---

[57] "Spatium Continuum componitur ex partibus qualibet a nobis determinabili minoribus". "Sunt aliqua in continuo infinite parva, seu infinities minora, quovis sensibili dato". "Punctum unum alio potest esse infinities minus".

[58] "Idem brevius: quodlibet ex his triangulis scalenis punctum est, cum eius basis et altitudo sint lineae infinite parvae. Ergo omnia simul, in unam rectam ex evoluta extensa conflatam disposita (tot enim sunt scalena, quot in ea chordae, seu puncta) non nisi lineam facient." Note that not only infinitely small lines are called points, but also infinitely small figures (cf. in particular *Plagula* 16₄ A7.4 pp. 304-331).

[59] "*Indivisibilia* definienda sunt infinite parva, seu quorum ratio ad quantitatem sensibilem (vel differentia) infinita est".

[60] It is interesting to note that a similar position, namely the idea that one can reach infinitesimals at an infinitieth stage of getting smaller and smaller, will be maintained by Bernoulli, and repeatedly criticized by the mature Leibniz, in their 1698 correspondence (GM III 529-554, see in particular pp. 539, 541 and 551).

[61] "Intelligatur abscissa AE dividi in partes aequales infinitas, quales sunt EF. FG. […] Eritque recta EF. vel GH. infinite parva, eademque erit infinitesima rectae AE. abscissae".

[62] "Demonstratio haec est: intelligatur triangulum figurae characteristicum esse MNE. cuius scilicet altitudo MN est infinitesima altitudinis figurae". Cf. A7.4 p. 432; 493; 494; 497; 663.

Infinitesimals are often employed as "unities":

> And we suppose EE infinitesimal, that is the unity for the construction = 1 [*Varia circa functiones...* 1673 A7.4 p.755][63]
> I'll grant you that, but I maintain that this 1 must be denoted in this way 1r. In fact it is not 1, or a certain number, but is something indivisible, which in this construction works as unit, that is to say the infinitesimal part of a certain line, divided in thought into infinite equal parts [*De paraboloeidum et hyperboloeidum quadratura II* 1673 A7.4 p. 607][64]

As in phase II, the notion of unassignable quantities does not seem particularly to worry Leibniz, but the new writings bear witness to a deepening in his analysis. On one hand, the vague notion of indistant parts disappears; on the other hand, we encounter the first mentions of imagination: a curve is divided "in thought" (*cogitatione*) into infinitesimal pieces.

Moreover, a hiatus is being created between what we perceive in physics, namely a continuous matter, made of assignable parts, and what is behind the appearance, that is motion and its unassignable *conatus*. The apparent and the deeper level seems to have two different mathematical counterparts.

## Phase IV (1675-6)

In phase IV, we witness the re-emergence of metaphysical explanations, which Leibniz had strongly criticized in phase I, but with a new awareness, for now metaphysics does not supplant mathematics, but will complement it. The idea is that mathematics provides a description of natural phenomena, which although approximate and idealized makes physical problems manageable and as predictable as possible. But mathematics does not give the true laws of physics, which are instead provided by metaphysics. As Leibniz recalls decades later:

> Mais quand je cherchay les derniers raisons du Mechanisme et les loix mêmes du mouvement, je fus tout surpris de voir qu'il etoit impossible de les trouver dans les Mathematiques, et qu'il falloit retourner à la Metaphysique. C'est ce qui me ramena aux Entelechies, et du materiel au formel. [*Leibniz to Remond*, 1714. GIII 606][65]

*Phase IV in physics*

Concerning the problem of the composition of physical bodies, Leibniz resorts to metaphysics in order to get the actuality he needs. Although physical bodies are infinitely divisible – namely, they are capable of all resolutions into parts – the potentiality of these parts is overcome by attributing to parts not only a motion but a soul, or mind:[66]

> And so one must examine rigorously whether there follows the perfect division of a liquid into metaphysical points, or only into mathematical points. For mathematical points could be called *indivisibiles* of Cavalieri, even if they are not metaphysical, i.e., are not minima. But if it could be shown that a liquid is divided to a greater or lesser extent, it would follow that a liquid is not resolved in[to] indivisibiles. But one could defend the composition of a liquid out of perfect points, even if it is never wholly resolved into them, since it is capable of all resolutions, and will cease when its cement – namely, mind and motion – ceases [*De Summa Rerum* 1676. A6.3. 474. Transl. Parkinson].[67]

---

[63] "Et posita EE infinitesima seu unitate constructionis = 1".

[64] "Ego fateor, sed aio illud 1 notandum esse hoc modo 1r. neque enim esse 1, seu numerum illum, sed esse indivisibile aliquod, quod in praesenti constructione unitatis personam sustinet, sive quod est infinitesima pars lineae cuiusdam, in partes aequales infinitas cogitatione divisae".

[65] "But when I looked for the ultimate reasons for mechanism, and even for the laws of motion, I was greatly surprised to see that they could not be found in mathematics but that I should have to return to metaphysics. This led me back to entelechies, and from the material to the formal" [Transl. Loemker].

[66] In this period Leibniz uses the Latin "mens" (mind) to mean every kind of mind-like substances. In later writings he distinguishes between minds (*mentes*), souls (*animae*) and "bare" monads. Monads are capable of perceptions; souls are capable of perception and consciousness; minds are capable of rational perceptions and self-consciousness.

And a single indivisible body – or ἄτομον – arises, of whatever magnitude, every time it has a soul [*Metaphysical notes* 1676. A6.3. 393].

In this way, the material parts of a physical body can be said to be unities in a proper sense – namely, indivisible actual wholes – and simultaneously infinite.

The situation is apparently similar to the one described in phases II and III: we have an infinitely divided material stuff, and "something" that keeps the parts together. The crucial difference is that while in phases II and III both the material parts and the "cement", that is their *conatus* of motion, were physical entities, in phase IV the bond is not physical but metaphysical, due to the fact that parts have minds, or souls.

So, while in phase III a geometrical representation of both the infinitely divisible assignable parts and their unassignable *conatus* was needed, because both belonged to physics, in phase IV a geometrical description is needed only for the infinitely divisible physical parts. As metaphysical notions, indeed, minds reside on another level.

In particular, they can be conceived as unassignable and unextended without any need to resort to a definition of "magnitude" different from "extension". In phase IV, in fact, being unextended is something more than a quantitative property: minds are unextended because, as any metaphysical entity, they are extraneous to the category of extension, and more generally to the category of quantity. As such, they are also unassignable.

If at the very beginning the separation between minds and their expression is not so clear,[68] as Leibniz himself admits in later writings, in his mature theory any conflation of metaphysical and material elements of a body disappears:

> Many years ago, when my philosophy was not yet sufficiently mature, I located souls in points […] But after further reflection, I discovered that in this way we are not only led into innumerable difficulties, but that there is also here a certain, as it were, category mistake (μετάβασιν εἰς ἄλλο γένος). Those things that pertain to extension should not be attributed to souls, and their unity and multitude should not be taken from the category of quantity, but rather from the category of substance.[69] [*Leibniz to Des Bosses* April 30 1709; G II, 372. Transl. Look-Rutherford. See Leibniz (2007)]

Leibniz's reference to a category mistake is significant because it points at the very heart of the problem: the conflation of substantial and quantitative analysis of physical bodies.

As we mentioned, in phase I Leibniz thought that he could avoid any difficulty of this kind by banishing metaphysics from physics: no substantial immaterial forms are needed to define physical objects, which are completely identified by their geometrical shape, and mechanical implications.[70] Physical objects are continuous, divisible into potentially infinite assignable extended parts, where the adjective *infinite* is taken in its quantitative meaning.

The main problem with this geometric/mechanical interpretation was the indefiniteness of the parts into which bodies can be analysed, but of which they cannot be "made", at least in a proper sense. In phases II and III Leibniz tried to solve the problem by hypothesizing actual ultimate elements, which held together by motion make a continuum, but have a proper unity. However, his move raised two new difficulties: the mathematical inconsistency of the model, which should bring together the actuality of the parts and their quantitative infinity, and the impossibility of distinguishing a unitary body from a simple aggregate of parts, or atoms, without resorting to some kind of form.

---

[67] "Adeoque examinandum rigorose, an perfecta divisio liquidi in puncta metaphysica, an vero tantum in puncta mathematica sequatur. Nam puncta mathematica possent appellari indivisibilia Cavaleriana, etsi non sint metaphysica seu minima. Quod si ostendi posset plus minusve dividi liquidum, sequetur non resolvi liquidum in indivisibilia. Posset tamen defendi liquidum componi ex punctis perfectis, etsi nunquam in illa prorsus resolvatur, vel ideo quia omnium est resolutionum capax, et caemento cessante, mente scilicet et motu, cessabit".

[68] Scholars' interpretations are discordant on this point. On the coexistence of matter and souls see for instance Garber 2009a. The thesis of a wholesale separation from the start is maintained in Mercer 2001.

[69] The relation between simple substances and souls is discussed in detail in Garber 2009a and 2016.

[70] In phase I, the only metaphysical notion Leibniz allowed was a supreme mind, which is the ultimate cause of motion, and hence of interaction between bodies.

Leibniz's minds are a first attempt to solve both problems from a new perspective. On one hand, in fact, minds give unity to the ultimate constituents of bodies, and keeping track of how the body has been created, make it distinguishable from other bodies with the same shape:

> We say that the effect involves its cause; that is, in such a way that whoever understands some effect perfectly will also arrive at the knowledge of its cause. For it is necessary that there is some connection between a complete cause and the effect. But on the other hand there is this obstacle: that different causes can produce an effect that is perfectly the same. For example, whether two parallelograms or two triangles are put together in the appropriate way (as is evident here) the same square will, as is clear, always be produced. Neither of these can be distinguished from the other in any way, not even by the wisest being. So, given a square of this kind, it will be in the power of no one – not even the wisest being – to discover its cause, since the problem is not determinate. The effect, therefore, seems not to involve its cause. So if we are certain, from some other source, that the effect does involve its cause, then it is necessary that the method of production must always be discernible in squares that have been produced. And so it is impossible that two squares of this kind should be perfectly similar; for they will consist of matter, but that matter will have a mind, and the mind will retain the effect of its former state [*De summa rerum* 1676. A6.3. 490-1. Transl. Parkinson p.51].[71]

On the other hand, as metaphysical objects, minds are non-quantitative notions, and the problem of their mathematical representation simply disappears. In fact, on one level we have the appearance of physical objects, that is the extended material bodies we perceive, and we can approximately represent using the quantitative objects of traditional geometry: a physical body will be idealized as a continuous, in(de)finitely divisible, geometric solid, every part of which is comparable with every other.

On a totally separate level we have the metaphysical foundation of the same physical object, which we do not directly perceive, and cannot represent using traditional geometry:

> Le corps est un aggrégé de substances, et n'est pas une substance à proprement parler. ll faut par conséquent que partout dans le corps il se trouve des substances indivisibles, ingénérables et incorruptibles, ayant quelque chose de répondant aux ames. Que toutes ces substances ont tousjours esté et seront tousjours unies à des corps organiques, diversement transformables. Que chacune de ces substances contient dans sa nature *legem continuationis seriei suarum operationum*, et tout ce qui luy est arrivé et arrivera. [*Leibniz à Arnauld 23 March 1690*, A2.2. 311-12] ("The body is an aggregate of substances, and is not a substance strictly speaking. It is therefore necessary that everywhere in the body there be substances that are indivisible, non-generatable and incorruptible, possessing something corresponding to souls. And that all these substances have always been and will always be attached to organic bodies, variously transformable. And that each of these substances contain in its nature *legem continuationis seriei suarum operationum*, and all that has happened to it and will happen to it.")

---

[71] "Dicimus effectum involvere causam suam; id est ita ut qui perfecte intelligat effectum aliquem, etiam ad causae eius cognitionem perveniat. Utique enim necessaria quaedam inter causam integram et effectum connexio est. Sed contra obstat, quod diversae causae producere possunt eundem perfecte effectum, exempli causa, sive duo parallelogramma, sive duo [triangula], uti oportet, et hic apparet, componantur, semper idem plane quadratum prodibit; quorum alterum ab altero nullo poterit discerni modo; ne a sapientissimo quidem; ita ut in nullius ne sapientissimi quidem potestate futurum sit, ex dato quadrato ejusmodi invenire eius causam, quoniam problema non est determinatum. Videtur ergo effectus non involvere suam causam. Quare si aliunde nobis certum sit, effectum involvere suam causam, necesse est in productis quadratis modum productionis semper discerni posse Adeoque impossibile esse, ut duo quadrata eiusmodi sint perfecte similia, quia ex materia constabunt, ea autem mentem habebit, et mens retinebit effectum status prioris". As noted by Garber in commenting this quotation: "This would seem to lead directly to the later view that individual substances, complete beings, require a soul as the seat of their complete individual concept: there must be form in bodies in order for them to contain the marks and traces of their previous states" [Garber 2009b, p. 72].

Summarizing, in Phase IV the transition occurs towards Leibniz's definitive theory of physical bodies, which are nothing but the phenomenal "expression" of the aggregation of infinite unextended metaphysical simple substances. Now, as we will show more in detail in our analysis of Phase V, the way in which physical bodies express metaphysical substances is only an aspect of a more general situation in which the whole physics is grounded on metaphysics, rather than on mathematics, as it was before.

*Phase IV in mathematics*
It is in mathematics that this separation has the most important consequences. In fact, once mathematics is no longer the foundation of physics, but at most provides an (approximate) description of it, it can be extended to include concepts which do not have any physical counterpart. In particular, lacking a physical model, these notions lack a geometrical representation too: since they cannot be shown in the phenomenal realm, they can only be imagined, or feigned.

In phases II and III, Leibniz's physical bodies were made of matter *and* motion: infinite assignable portions of a continuous material stuff held together and made actual by their unassignable endeavors of motion. A geometrical representation that would account consistently for both aspects turned out to be problematic, having to make the actuality of the unassignables coexist with the infinite divisibility of the continuum in assignable parts.

Gradually, Leibniz abandoned the attempt to keep them together: physical bodies have a natural representation as continuous geometrical extensions, while unassignable infinitely small quantities become imaginary notions, or fictions, and we will conclude the analysis of Phase IV following the path of Leibniz's infinitely small quantities from traditional geometry towards fictionality.

The crucial step of this path is the awareness of the consequences of unassignability, namely the non-Archimedean nature of infinitely small quantities. In particular, being in no ratio to finite magnitudes, when compared with them infinitely small quantities can be ignored, without leading to any contradiction.[72]

We argue that the syncategorematic reading represents a step in Leibniz's path toward his mature conception of infinitesimals, but not the final result. This is evident from the reading of the notes he wrote in the years 1675/76, and in particular of *De Quadratura Arithmetica* (*DQA*), Leibniz's longest mathematical work.

In DQA three methods for the quadrature of a curve are employed: the traditional indirect Archimedean method of exhaustion, Cavalieri's method of indivisibles, and Leibniz's infinitesimal techniques. In his seminal study of *DQA* Eberhard Knobloch considered as substantially equivalent Leibniz's and Cavalieri's method, and argued that indirect Archimedean methods are furnished so as to justify these speedy but risky techniques. Leibniz himself warns the reader concerning the risks inherent in the method of indivisibles:

---

[72]   The only requirement is that in deleting, and more generally in manipulating infinitely small quantities, one follows certain well-defined rules; see note 22.

> J'ai pris plaisir à poursuivre ce raisonnement car il offre un exemple des précautions à prendre lorsqu'on raisonne sur l'infini dans la méthode des indivisibles [*DQA* 1676 Prop. XXII *scholium* A 7.6 583, translated by Parmentier p. 177].[73]

("I took pleasure in pursuing this reasoning because it offers an example of the precautions to be taken when reasoning about infinity in the method of indivisibles.")

We argue that infinitesimal methods are rigorous, while according to Leibniz, Archimedean methods must be handled with care. In fact, if we have to apply an exhaustion-like procedure to an infinite area, we must be careful, because some absurd results could appear.

With the introduction of *infinita terminata*, and the rules for computing with them (see footnote 83 below), Leibniz provides a tool for escaping such difficulties,

> et montre qu'on ne peut pas toujours s'empresser d'affirmer que la totalité d'un espace infini possède une certaine propriété, sous prétexte que celle-ci est toujours vérifiée par un nombre fini de parties [*DQA* 1676 Prop. XXII *scholium* A 7.6 583, translated by Parmentier p. 177].[74]

("and shows that one cannot always hasten to affirm that the totality of an infinite space possesses a certain property, under the pretext that this property is always verified by a finite number of parts.")

The DQA is sometimes used for supporting the claim that Leibniz's mature infinitesimals are inherently syncategorematic, so that the occurrences of the term "infinitesimal" are are merely signs for traditional exhaustion procedures, or else foreshadow "epsilontic" ones.[75] However, an attentive reading of DQA shows the opposite. In fact, it can be divided in two parts, containing respectively a traditional syncategorematic and a new infinitesimal approach.[76]

The first part of DQA (Prop I-X) does not contain anything new concerning infinitely small and infinite quantities. All one needs in order to understand propositions I-X is the syncategorematic stuff typical of Greek geometry: we take a finite quantity and then we divide (or we multiply) it indefinitely. We can go on indefinitely due to the properties of the continuum : at every step we have a finite quantity, which at every step is smaller (bigger), without end. Both the notion of the asymptote (a line that can be made bigger and bigger without meeting the curve) and that of negligible error (a quantity that can be made lesser and lesser without becoming zero) belong to this traditional approach. Cavalieri's method of indivisibles is thought by Leibniz to fit into this scheme: a curvilinear area is filled with rectangles with infinitely small basis; since one can show that, diminishing the basis, the remaining part can be made lesser than any given quantity, the equivalence between the "sum" of infinite rectangles and the curvilinear area is demonstrated (in fact this is closer to Torricelli's method than Cavalieri's).

---

[73] "Libenter hanc contemplationem persecutus sum, quia specimen exhibet cautionis circa ratiocinia de infinito, et methodum indivisibilium"

[74] "ostenditque non semper ex partium finitarum perpetuo abscissarum proprietate quadam ad totius infiniti spatii proprietatem posse prosiliri"

[75] The idea of "syncategorematic infinitesimals" – a locution that never appears in Leibniz known writings but has been introduced in Arthur 2013 in order to explain Leibniz's mature conception of infinitesimals – comes from the idea of a syncategorematic infinite: just as a syncategorematic infinite refers to the procedure of indefinitely increasing a finite quantity, which remains in any case finite (see Section 2.2), syncategorematic infinitesimals refer to the procedure of indefinitely decreasing a finite quantity, which remains in any case finite. In other words, infinitesimals are a shorthand for Archimedean exhaustion-style procedures. It is clear that the definition itself of syncategorematic infinitesimals depends on the assumption that for Leibniz the counterpart of infinitesimals are unbounded infinities. Based on Leibniz's own texts, we reject this assumption. Similarly, we reject the idea that Leibniz's infinitesimals foreshadow the current traditional theory of limits, usually referred as "epsilontic". On the epsilontic interpretation of Leibniz's infinitesimal calculus and its history see in particular Breger 2008, pp. 187-189: "The connection between infinitesimals and what we now call epsilontics was *obvious enough* for 17th-century mathematicians" (p. 187).

[76] In fact, Leibniz introduces two novelties in DQA: (1) he introduces some "fictions" i.e., lines which are infinite and bounded. As we read in Prop. XI: "dans la mesure où elle s'appuie sur une fiction et aboutit sans difficulté dès qu'a été admise une ligne à la fois bornée (terminata) et infinie [DQA, Prop XI, p.99]. But also (2) he generalizes the method of indivisibles (see prop. VI and Parmentier's footnotes, especially n. 2 p. 63). And in fact he calls Cavalieri's method "la methode des indivisibles ordinaire"). Here the novelty consists in the introduction of what Parmentier calls "rectangles intermédiaires, dont les arêtes ne correspondent pas aux points de la courbe" (DQA note 2 p.63).

In most cases, the traditional method works perfectly, but there are instances where, as Leibniz explicitly says, it is not safe. In particular, it becomes "slippery" (lubrica) when one has to consider an infinite area, for instance the one contained between a curve and its asymptote. In this case, one can reach a contradiction, and in Proposition XXII he shows a typical case: an infinite area is proved to be equal to the same area minus a finite part:[77]

> On aboutira de la même manière, pour les autres hyperboloïdes, à un tout égal au tiers ou au quart d'une de ses parties. Voici une élégante manière de montrer combien sont périlleux les raisonnements sur l'infini lorsqu'ils ne sont pas pour ainsi dire guidés par le fil d'une démonstration [*DQA* 1676 Prop. XXII *scholium* A 7.6 583, translated by Parmentier pp. 177-179].[78]

("We will arrive in the same fashion, for the other hyperboloids, at a whole equal to a third or a quarter of one of its parts. Here is an elegant way of showing how perilous reasoning about infinity is when it is not, so to speak, guided by the thread of a demonstration.")

It is worth noting that the contradiction follows from having assumed that the quantities involved in the proof are always finite, at any step of the procedure. But this is exactly the crucial property of syncategorematic procedures, and of the unbounded infinities they involve.[79]

By contrast, starting from prop. XI, Leibniz lays the foundation of a new method, distinct from the traditional one, and which clearly has the purpose of avoiding ambiguity and contradictions related to the use of unbounded infinities:[80]

> Bien plus, à supposer que les questions dont je traite soient susceptible de démonstrations directes, je n'hésiterais pas à avancer qu'on n'en pourrait donner aucune qui n'admette ces quantités fictives, infiniment petites ou infinies, reportez-vous plus haut au scholie de la proposition VII [81]

---

[77] "Eodem modo in aliis colligetur totum partis suae trientem aut quartam partem esse. Eleganti argumento quam lubrica sit ratiocinatio circa infinita, nisi demonstrationis filo regatur" (Translated by Parmentier). Historians often claim that it is Leibniz's method with infinitesimals that is not safe and therefore requires justification by means of an exhaustion like procedure (on this interpretative line see Rabouin-Arthur 2020). In fact here Leibniz is claiming on the contrary that only his direct method is safe, which does not involve unbounded infinities. There is currently a lively debate concerning the Leibnizian calculus; see e.g., Bair et al. (2021) and references cited therein, as well as Bair et al. (2022), Katz et al. (2022), Archibald et al. (2022), Bair et al. (2023), Katz et al. (2023).

[78] "Eodem modo in aliis colligetur totum partis suae trientem aut quartam partem esse. Eleganti argumento quam lubrica sit ratiocinatio circa infinita, nisi demonstrationis filo regatur". Here Leibniz is clearly referring to arguments using infinite wholes, rather than arguments using infinitesimals.

[79] Although this demonstration has been used to argue the inherently contradictory nature of Leibniz's *infinita terminata* (see for instance Arthur 2013, pp. 556-558, where an early version of the proof is used) it actually proves the contradictory nature of *infinita interminata*, rather than of *infinita terminata*.

[80] Significantly, Knobloch admits that contradictions are due to the use of, specifically, *un*bounded infinity: "Leibniz underlines that the method of infinite quantities is dangerous. We have to apply the thread of an indirect proof if we have in mind to rely on it. Why? Because the infinite is characterized in a privative way. In reality, there is no last point, no end, we have only relied on a fiction in order to admit it in mathematics." Knobloch 1994 (p. 270). But the fiction is not the slippery unbounded infinite: it is the reliable bounded infinite.

[81] Although they are not defined until Prop. XI, fictive infinite quantities are mentioned in the *Scholium* to Prop.VII. This fact seems to have led some readers to conflate the two methods. But in Prop VII Leibniz does not use any fiction: he merely notes fact that he has been forced to use a *reductio ad absurdum* because he does not know any direct method. In fact, he fears that no direct solutions are possible, unless one agrees to use fictive quantities, namely infinitesimals and *infinita terminata*: "j'ai même des raisons de craindre qu'on ne puisse pas le faire d'une manière

("Furthermore, assuming that the questions I am treating are susceptible of direct proofs, I would not hesitate to assert that one could not give any [proofs] which do not exploit these fictitious quantities, infinitely small or infinite, see above at Scholium of Proposition VII.")

> [...] Le présent petit traité en fournira un exemple propre à convaincre de la fécondité de la méthode; ce que j'ai dit dans le scolie de la proposition 22 fournira d'autre part un exemple particulier de raisonnement sûr, là où la méthode de Cavalieri prise au sens strict n'est pas fiable [*DQA* 1676 Prop. XXIII *scholium* A 7.6 585-586, translated by Parmentier][82].

("This little treatise will provide a suitable example to convince of the fecundity of the method; what I said in the Scolium of Proposition 22 will, on the other hand, furnish a particular example of reliable reasoning, where Cavalieri's method taken in the strict sense is not reliable.")

Here Leibniz furnishes an example where the traditional method is not reliable, whereas the introduction of Leibniz's fictions furnishes a more reliable method.

The main ingredients of Leibniz's method are bounded infinities (*infinita terminata*). Here the infinitely small quantities ($\mu$) are fractions whose denominators are infinitely large but bounded quantities ($\mu$)$\lambda$, both called fictions. It does not matter – Leibniz says – if these quantities are *in rerum natura*. Of course, for using them safely, one needs to provide the appropriate rules, indeed provided by Leibniz.[83]

These quantities (bounded infinities and infinitesimals) in the DQA are treated as different from what Leibniz calls indivisibles. In particular, Leibniz says that:

> Au demeurant il semblera à certains paradoxal de parler de Lignes infinies qui soient cependant bornées, je dois donc faire remarquer qu'il existe une grande différence entre indivisible et infiniment petit, et de la même manière une grande différence entre infini et non borné [*DQA* 1676 **Prop. XI**][84]

("At any rate, it will seem paradoxical, to some, to speak of infinite Lines which are nevertheless bounded, I must therefore point out that there is a great difference between indivisible and infinitely small, and in the same way a great difference between [bounded] infinite and unbounded.")

Here Leibniz emphasizes the difference between indivisibles and infinitesimals, which supports our interpretation. Leibniz also wrote:

---

naturelle et sans faire intervenir des quantités fictives, je veux dire infinies ou infiniment petites" (imo ratione habeo, cur verear ut id fieri possit per naturam rerum sine quantitatibus fictitiis, infinitis scilicet vel infinite parvis assumptis *DQA* VII scholium, translated by Parmentier).

[82] "Imo si quidem possibile est directas de his rebus exhiberi demonstrationes, ausim asserere, non posse eas dari, nisi his quantitatibus fictitiis, infinite parvis, aut infinitis, admissis. Adde supra prop. 7. schol. [...] Cujus specimen totus hic libellus erit si quis methodi fructum quaerit; securitatis autem excemplum peculiare dabunt quae diximus *schol. ad prop 22* ubi Cavalerii methodus crude sumta infida est". What may have misled some readers is Leibniz's reference to the method of indivisibles. They take this to be a reference to infinitesimals, whereas we argue that Leibniz is referring to exhaustion-type arguments.

[83] In particular, having constructed the infinite rectangle with basis ($\mu$) and length ($\mu$)$\lambda$ Leibniz proves that, depending on the properties of the curve, it can be: (a) finite, (b) infinite or (c) infinitely small. The proof of a is made extending to the infinite rectangle some properties which hold for finite areas. The proofs of (b) and (c) are made using the rectangle xy where x is the "first ordinate" and y is the "first abscissa". Summarizing, Leibniz shows that if $x : l = m : (\mu)\lambda$ with ($\mu$)$\lambda$ infinitum terminatum, then x is infinitely small. if ($\mu$) : l = m : y with ($\mu$) infinitely small, then y is infinite (terminatum).

[84] Note that the contrast is a chiastic one: indivisibles and unbounded go together, and infinitesimal and infinite (i.e., bounded infinite) go together, as is clear from the subsequent discussion.

> La Géométrie des indivisibles est trompeuse si on n'y adjoint quelques précisions au sujet des infiniment petits; utiliser des points véritablement indivisibles comporte un certain risque et il faut au contraire employer des lignes qui, pour être infiniment petites, n'en sont pas moins des lignes et sont, par conséquent, divisibles [*DQA* 1676 **Prop. XI**].

("The Geometry of the Indivisibles is misleading if one does not add to it some clarifications concerning the infinitely small; using truly indivisible points carries a certain risk and on the contrary it is necessary to use lines which, fwhile being infinitely small, are nonetheless lines and are, therefore, divisible.")

This clearly establishes that when Leibniz is talking about indivisibles, he is not talking about infinitesimals.

In sum, in *De Quadratura Arithmetica* Leibniz presents the infinitesimal method as more reliable than the Archimedean one, because the latter sometimes relies on infinite wholes which are contradictory. In better-known passages, however, he justifies his infinitesimal method by an appeal to Archimedean exhaustion. It is possible that the erroneous idea that in *De Quadratura Arithmetica* Leibniz speaks of unreliability of his method, rather than of other methods involving unbounded infinities, comes from these later passages.

In phase IV, one is far from the unproblematic approach of phases II and III, when Leibniz was convinced that, although unassignable, infinitely small quantities could simply be handled via the method of exhaustion. Now he feels the need to justify their legitimacy.

As sketched in our analysis of the *De Quadratura Arithmetica,* Leibniz's strategy consists of two steps. On one hand he shields his unassignable infinitely small quantities from the accusation of not complying with the part-whole axiom, and in this way he renders them consistent from a mathematical point of view. On the other hand, he introduces for them a peculiar *status*, that of "fictions", which holds together their mathematical consistency and their lacking of any physical counterpart.

Concerning the part-whole axiom, Leibniz specifies that only unbounded infinities (*infinita interminata*) contradict it, and for this reason are absolutely impossible.[85] By contrast, bounded infinities imply no contradiction, and for this reason are only accidentally impossible.[86]

Concerning the *status* of his unassignable infinitely small quantities, Leibniz assimilates them to other physically problematic but useful, well defined and well-known notions such as imaginary numbers, which appear in mathematics but cannot be exhibited or geometrically constructed, as we have shown at the end of phase III and as reaffirmed by Leibniz in phase IV in the following terms:

> I call *imaginary constructions* the constructions that we can feign (fingere), but we can show (exhibere) only with your mind. They are nevertheless possible […]: we can imagine for instance a quantity that grows – or is shortened or is transformed or is moved – following a given law, to which apparently no "ratio" can be associated. And it can happen that one can find a physical or even a geometrical realisation (*effectio*) of a construction that at the beginning seemed imaginary [*De descriptione geometrica curvarum transscendentium* 1675 A7.7 516].[87]

---

[85] We know that already in 1672 Leibniz declared infinite wholes to be contradictory: "Hence it follows either that in the infinite the whole is not greater than the part, which is the opinion of Galileo and Gregory of St. Vincent, and which I cannot accept; or that infinity itself is nothing, i.e., that it is not one and not a whole". [*Notes on Galileo's Two New Sciences* 1672. A VI. 3, p. 168; translated in Arthur 2001, p. 9].

[86] "Therefore if the essence of a thing can be conceived, provided that it is conceived clearly and distinctly (e.g., a species of animal with an uneven number of feet, also a species of immortal beast), then it must already be held to be possible, and its contrary will not be necessary, even if its existence may be contrary to the harmony of things and the existence of the divine intelligence, and consequently it never will actually exist, but it will remain per accidens impossible. Hence all those who call impossible (absolutely, i.e., per se) whatever neither was nor is nor will be are mistaken". (*Confessio Philosophi* A VI 3, 128, transl. Sleigh).

> One must discuss if √−1 is nothing or if it contains a truth. In fact, even if it cannot be done, it can be thought, not in itself but using characters and analogy (I call this way of thinking a blind one). In fact, as incommensurables are such that their powers are commensurable, in the same way imaginary are such that their powers are real, and impossibles are such that their squares are possibles, like √−1, whose square is -1, even if we assume that in nature there is no quantity corresponding to that: in fact, it is enough that its character is useful, because together with other ones it can express something real [*Imaginariae usus ad comparationem Circuli et Hyperbolae* 1675 A7.7 561].[88]

Like imaginary numbers, unassignable infinitely small quantities are intrinsically not geometrically representable: how could one represent a non-Archimedean quantity in a continuous framework such as that of Greek geometry? Furthermore, they cannot be viewed as the idealisation of some physical objects, which would not respect the harmony of things. Summarizing: since unassignable infinitely small quantities do not contradict the part-whole axiom, they are possible (more precisely, not absolutely impossible) notions, but since they cannot be exhibited – either in geometry or in physics – they are fictions. The same position will be reiterated by Leibniz in subsequent years:

> And some extractions of roots are such that roots are surd and they do not exist in natura rerum, and we call them imaginary, and the problem is impossible, as when our analysis shows that the requested point must be exhibited by the intersection of a specific circle and a specific straight line, in which case it may happen that this circle by no means reaches this line, and then the intersection is imaginary […] There is a big difference between imaginary quantities, or those impossible by accident, and absolutely impossible ones, which involve a contradiction: e.g., when it is found that solving a problem requires that 3 be equal to 4, which is absurd. But imaginary quantities, or [quantities] impossible by accident, namely quantities that cannot be exhibited for lack of a sufficient condition, which is required for having an intersection, can be compared with infinite and infinitely small quantities, which are generated in the same way […] And it is true that calculus necessary leads to them, and people who are not sufficiently expert in such matters get entangled and think they have reached an absurdity (*absurdum*). Experts know instead that this apparent impossibility [apparentem illam impossibilitatem] only means that a parallel line is traced instead of a straight line making the required angle, and this parallelism is the required angle, or quasi-angle [*Elementa nova matheseos universalis* 1683 A 6.4 520-521][89]

---

[87] "*Constructiones imaginarias* voco, quas fingere possumus, exhibere non nisi mente non possumus. Sunt tamen possibiles […] tales sunt si imaginemur magnitudinem aliquam certa lege crescere aut contrahi aut transformari, aut etiam moveri, cujus reapse praescribendae nulla appareat ratio. Et fieri poterit aliquando ut constructionis quae initio imaginaria apparuit tandem physica aut etiam Geometrica effectio reperiatur". The passage was discussed also in Katz et al. (2023).

[88] "An autem quantitas √−1. sit omnino nihil, an vero contineat nescio quid discutiendum diligentius. Etsi enim non possit effici, potest tamen quodammodo intelligi, non seipsa sed ope characteris et analogiae, exemplo illius cogitationis quam caecam appello. Et vero quemadmodum sunt incommensurabiles quae sunt potentia commensurabiles, ita sunt imaginariae quarum potentiae sunt reales; seu impossibiles quarum quadrata sunt possibilia, ut √−1, cujus quadratum −1 etsi autem ponatur nihil omnino esse in natura tali quantitati respondens; sufficit tamen ejus characterem esse utilem, quia cum aliis junctus realia exprimit".

[89] "Et quaedam extractiones tales sunt, ut radices illae surdae nec in natura rerum extent, tunc dicuntur imaginariae, et problema est impossibile, ut cum analysis ostendit punctum quaesitum debere exhiberi per intersectionem certi circuli et certae rectae, ubi fieri potest ut ille circulus ad illam rectam nullo modo perveniat, et tunc intersectio erit imaginaria […]. Multum autem interest inter quantitates imaginarias, seu impossibiles per accidens, et impossibiles absolute quae involvunt contradictionem, ut cum invenitur ad problema solvendum opus esse, ut fit 3 aequ. 4 quod est absurdum. Imaginariae vero seu per accidens impossibiles, quae scilicet non possunt exhiberi ob defectum sufficientis constitutionis ad intersectionem necessariae, possunt comparari cum Quantitatibus infinitis et infinite parvis, quae eodem modo oriuntur […] quod adeo verum est ut saepe calculus ad eas necessario ducat, ubi harum rerum nondum satis periti mire torquentur et in absurdum se incidisse putant. Intelligentes vero sciunt apparentem illam

The idea is that traditionally Greek geometry was considered the natural idealisation of physics, and Greek geometry did coincide with the whole mathematics.[90] Subsequently, mathematical notions were found which did not belong to Greek geometry – roughly speaking, Greek geometry became a subset of a larger mathematics – but the idea survived that the idealisation of physics has to be Greek geometry. The elements belonging to the complement of Greek geometry, having no direct physical counterpart, are then accepted, but they are regarded as fictions, or imaginary notions.

## Phase V (1677-...)

Phase V leads Leibniz to his ultimate theory of the composition of bodies, a theory in which the separation among mathematics, physics and metaphysics is complete. In particular, Leibniz eliminates any residual ambiguity by providing a metaphysical foundation to any feature of a physical body. Not only the actuality and the unity of the parts, but their physical extension too is now explained in term of metaphysical notions.

Basically, Leibniz's idea is to extend to physical composite substances his analysis of simple metaphysical substances. Simple metaphysical substances – which the mature Leibniz also calls monads – can be analysed into an active force (*vis activa*), or substantial form, and a passive force, or prime matter.[91] In an analogous way, in phase V Leibniz expresses both extended matter and substantial form of physical bodies in terms of forces: a primitive passive force, corresponding to matter – the "*primary matter* in the schools, if correctly interpreted"–,[92] and a primitive active force, corresponding to substantial form. As for Aristotle form and matter, together, constitute a substance, for Leibniz primitive forces, active and passive, do the same:

> Primitive active force, which Aristotle calls first entelechy and one commonly calls the form of a substance, is another natural principle which, together with matter or [primitive] passive force, completes a corporeal substance. This substance, of course, is one *per se*, and not a mere aggregate of many substances, for there is a great difference between an animal, for example, and a flock [*On Body and Force* (1702), AG 252].[93]

As Aristotelian primary matter and form, primitive forces are causes in a wide sense, for they are responsible for the whole being of a physical substance, and consequently for its behaviour, and its observable properties:

---

impossibilitatem tantum significare, ut loco rectae angulum quaesitum facientis ducatur parallela; hunc parallelismum esse angulum illum seu quasi angulum quaesitum."

[90] Still in Leibniz's days, algebraic notions were viewed as a sort of symbolic "translation" of geometrical quantities, and a "rigorous" proof meant a geometrical proof. On this point see Breger 2008, §3. Cf. Parmentier 1995, pp. 26-27.

[91] While prime matter, as pure passivity, cannot be further specified, Leibniz describes the active power, or *entelecheia*, of a simple substance in terms of a law, the *law of the series*, which expresses the inner states of the substance itself (see GP III, 66; 356; GP IV, 512). To explain the way in which the law of the series rules the sequence of the inner states, Leibniz uses the analogy of a mathematical sequence of numbers: "L'essence des substances consiste dans la force primitive d'agir, ou dans la loy de la suite des changemens, comme la nature de la series dans les nombres (*Excerpta ex notis meis marginalibus ad Fucherii responsionem primam in Malebranchium critica*, 1676, A6.3. 326). On the relation between the law of the series and simple substances in general see for instance Whipple 2010 and the bibliography quoted there.

[92] *Specimen Dynamicum* [GM VI, 237].

[93] "Vis activa primitiva quae Aristoteli dicitur ἐντελέχεια ἡ πρώτη, vulgo forma substantiae, est alterum naturale principium quod cum materia seu vi passiva [primitiva] substantiam corpoream absolvit, quae scilicet unum per se est, non nudum aggregatum plurium substantiarum, multum enim interest verbi gratia inter animal et gregem" (GP IV, 395).

> Composite substance does not formally consist in monads and their subordination, for then it would be a mere aggregate or a being per accidens. Rather, it consists in primitive active and passive force, from which arise the qualities and the actions and passions of the composite which are discovered by the senses, if they are assumed to be more than phenomena. [*Leibniz to Des Bosses* (1716), see Leibniz (2007), AG 203-4][94]

In particular, mechanical properties of physical substances – namely their internal coherence and their reciprocal interactions, and therefore the local motions behind them – are due to peculiar properties of these primitive forces, which Leibniz calls derivative forces:

> By derivative force, namely, that by which bodies actually act on one another or are acted upon by one another, I understand, in this context, only that which is connected to motion (local motion, of course), and which, in turn, tends further to produce local motion. For we acknowledge that all other material phenomena can be explained by local motion [*Specimen Dynamicum* (1695), AG 120].[95]

As for metaphysical notions in general, we have direct access neither to primitive forces, nor to derivative ones: at most, they can be partially understood by means of analogical speech. In particular, if we read primitive forces, active and passive, as substances, then derivative forces, active and passive, must be read as their faculties (*facultés*), or ways of being (*façons d'être*).[96] If otherwise we read the primitive force in analogy with the law of the series (of the states of the substance; see note 91), then derivative forces can be related to the single terms of the series:

> Derivative force is itself the present state when it tends toward or preinvolves a following state, as every present is pregnant with the future. But that which persists, insofar as it involves all cases, contains primitive force, so that primitive force is, as it were [*velut*], the law of the series, while derivative force is, as it were, a determination which designates some term in the series [*Leibniz to De Volder* (1704). Transl. Loemker].[97]

At this point Leibniz differentiates his derivative active force from Scholastics' faculties, which are traditionally interpreted as purely passive powers (*dynameis*). Derivative active force, instead, is endowed with an active tendency, which naturally tends towards its fulfillment:

> Active force differs from the mere power familiar to the Schools, for the active power or faculty of the Scholastics is nothing but a close possibility of acting, which needs an external excitation or a stimulus, as it were, to be transferred into action. Active force, in contrast, contains a certain act or entelechy and is thus midway between the faculty of acting and the act itself and involves a *conatus*. It is thus carried into action by itself and needs no help but only the removal of an

---

[94] "Substantia composita non consistit formaliter in monadibus et earum subordinatione, ita enim merum foret aggregatum, seu ens per accidens, sed consistit in vi activa, et passiva primitiva, ex quibus oriuntur qualitates et actiones passionesque compositi, quae sensibus deprehenduntur, si plus quam phaenomena esse ponantur" (GP II, 517-518).

[95] "Vim ergo derivativam, qua scilicet corpora actu in se invicem agunt aut a se invicem patiuntur, hoc loco non aliam intelligimus, quam quae motui (locali scilicet) cohaeret, et vicissim ad motum localem porro producendum tendit. Nam per motum localem caetera phaenomena materialia explicari posse agnoscimus" (GM VI, 237).

[96] "Les Puissances primitives constituent les substances mêmes, et les puissances derivatives, ou si vous voulés, les facultés, ne sont que des façons d'estre, qu'il faut deriver des substances (*Nouveaux Essais*, IV.iii.6. A VI, 6, 379).

[97] "Vis autem derivativa est ipse status praesens dum tendit ad sequentem seu sequentem prae-involvit, uti omne praesens gravidum est futuro. Sed ipsum persistens, quatenus involvit casus omnes, primitivam vim habet, ut vis primitiva sit velut lex seriei, vis derivativa velut determinatio quae terminum aliquem in serie designat" (GP II, 262).

impediment. [*De primae philosophiae emendatione, et de notione substantiae*, 1694. Transl. Loemker][98]

The explanation is resumed at the beginning of the *Specimen Dynamicum*, where Leibniz insists on the notion of *conatus*, now called *nisus* and defined in a Scholastic way as a disposition, or appetite, which completes the force:

> Elsewhere we urged that in corporeal things there is something over and above extension, in fact, something prior to extension, namely, that force of nature implanted (*inditam*) everywhere by the Creator. This force does not consist in a simple faculty, with which the schools have been content, but is further endowed with *conatus* or *nisus*, attaining its full effect unless it is impeded by a contrary conatus [*Specimen Dynamicus* 1695, AG 118]. [99]

Summarizing, in phase V Leibniz reads extended physical bodies as the mere phenomenal "expression" of the aggregation of infinite unextended metaphysical simple substances. Here a simple substance simply means the union of an active force, or substantial form, and a passive force, or prime matter. On a purely metaphysical level, then, physical bodies too can be analyzed in terms of primitive forces, active and passive. In their turn, primitive forces are meant to act in virtue of certain faculties, which Leibniz calls derivative forces, active and passive. Finally, derivative forces "produce" all the observable physical effects. For instance, all manifestations of local motion come from derivative active forces, via their *conatus*, or *nisus*.

The actuality of derivative forces manifests itself in observable but phenomenal local motions:

> For, strictly speaking, motion (and likewise time) never really exists, since the whole never exists, inasmuch as it lacks coexistent parts. And furthermore, there is nothing real in motion but a momentary something which must consist in a force striving *[nitente]* toward change. Whatever there is in corporeal nature over and above the object of geometry or extension reduces to this [**AG 119**].

In a similar way, the unity and self-sufficiency of metaphysical substances, namely of primitive forces, has its phenomenal expression in the physical extension of bodies, every part of which has a proper individuality.[100] In other words, due to its monadic origin, physical matter is not a continuous stuff, like geometrical extension is, but a contiguous one, every part being characterized and individuated by a different motion:

> As for the first point, it follows from the very fact that a mathematical body cannot be analyzed into primary constituents that it is also not real but something mental and designates nothing but the possibility of parts, not something actual. A mathematical line, namely, is in this respect like arithmetical unity; in both cases the parts are only possible and completely indefinite. A line is no more an aggregate of the lines into which it can be cut than unity is the aggregate of the fractions into which it can be split up. […]
> But in real things, that is, bodies, the parts are not indefinite – as they are in space, which is a mental thing – but actually specified in a fixed way according to the divisions and subdivisions

---

[98] "Differt enim vis activa a potentia nuda vulgo scholis cognita, quod potentia activa Scholasticorum, seu facultas, nihil aliud est quam propinqua agendi possibilitas, quae tamen aliena excitatione et velut stimulo indiget, ut in actu transferatur. Sed vis activa actum quendam sive ἐντελέχειαν continet, atque inter facultatem agendi actionemque ipsam media est, et conatun involvit; atque ita per se ipsam in operationem fertur; nec auxiliis indiget, sed sola sublatione impedimenti" (GP IV 469).

[99] "In rebus corporeis esse aliquid praeter extensionem, imo extensione prius, alibi admonuimus, nempe ipsam vim naturae ubique ab Authore inditam, quae non in simplici facultate consistit, qua Scholae contentae fuisse videntur, sed praeterea conatu sive nisu instruitur, effectum pleno habituro, nisi contrario conatu impediatur" (GM VI, 235).

[100] In Leibniz's words "every part of matter is actually divided into other parts" [*Leibniz to Arnauld*, November 28-December 8, 1686. A2.2. 122. Transl. Ariew-Garber]; "In truth, matter is not a continuous thing but a discrete one, actually divided at infinity" [*Leibniz to De Volder*, October 11 1705; G II, 278. (A2.4. 303600). Transl. Loemker].

which nature actually introduces through the varieties of motion. And granted that these divisions proceed to infinity, they are nevertheless the result of fixed primary constituents or real unities, though infinite in number. Accurately speaking, however, matter is not composed of these constitutive unities but results from them. [*Leibniz to De Volder* June 30 1704; G II, 268. Transl. Loemker][101]

In truth, matter is not a continuous thing but a discrete one, actually divided *ad infinitum* [*Leibniz to De Volder*, October 11, 1705; G II, 278].[102]

And certainly every extended thing is divisible, in such a way that parts can be assigned in it, but in matter they are actually assigned, while in extension they are only potential, as in number. [*Ad Schedam Amaxariam*].[103]

The crucial novelty, which solves the last remaining difficulties in Phase IV, is that now all the ultimate "components" of physical bodies – namely, primitive and derivative forces, as well as their tendencies – are metaphysical entities. And as metaphysical entities, the question of their geometrical representation simply disappears:

Therefore, I concluded from this that, because we cannot derive all truths concerning corporeal things from logical and geometrical axioms alone, that is, from large and small, whole and part, shape and position, and because we must appeal to other axioms pertaining to cause and effect, action and passion, in terms of which we can explain the order of things, we must admit something metaphysical, something perceptible by the mind alone over and above that which is purely mathematical and subject to the imagination, and we must add to material mass [massa] a certain superior and, so to speak, formal principle. Whether we call this principle form or entelechy or force does not matter, as long as we remember that it can only be explained through the notion of forces [*Specimen Dynamicum* (1695), AG 125].[104]

---

[101] "Quod primum attinet, eo ipso quod corpus mathematicum non potest resolvi in prima constitutiva, id utique non esse reale colligatur, sed mentale quiddam nec aliud designans quam possibilitatem partium, non aliquid actuale. Nempe linea mathematica se habet ut unitas arithmetica, et utrobique partes non sunt nisi possibiles et prorsus indefinitae; et non magis linea est aggregatum linearum in quas secari potest, quam unitas est aggregatum fractionum in quas potest discerpi […] At in realibus, nempe corporibus, partes non sunt indefinitae (ut in spatio, re mentali), sed actu assignatae certo modo, prout natura divisiones et subdivisiones actu secundum motuum varietates instituit, et licet eae divisiones procedant in infinitum, non ideo tamen minus omnia resultant ex certis primis constitutivis seu unitatibus realibus, sed numero infinitis. Accurate autem loquendo materia non componitur ex unitatibus constitutivis, sed ex iis resultat".

[102] "Revera materia non continuum sed discretum est actu in infinitum divisum".

[103] "Et sane omne extensum divisibile est, ut partes in eo assignari possint, actuque sint assignatae in materia, potentiales vero sint in ipsa extensione, ut in numero". A similar statement is reiterated a few pages later: "Est enim ut saepe dixi dispositio compossibilium phaenomenorum, Geometriaque est possibilitatum ipsarum tractatio. Et ideo extensio ipsa possibilis sive continua non habet partes determinatas, non magis quam Unitas numerica; quod secus est de physicis, seu existentibus ubi divisiones sunt actu factae. Quemadmodum non datur minima fractio in numeris, ita nec minimum in Geometria, punctumque est extremum non pars. Sed in rebus ipsis dantur minima simplices nempe substantiae. (*Ad Schedam Hamaxariam*, ad. 16, in LH IV, 3, 5c, Bl. 1-2. Transcription by Massimo Mugnai and Heinrich Schepers).

[104] "Hinc igitur, praeter pure mathematica et imaginatione subjecta, collegi quaedam metaphysica solaque mente perceptibilia esse admittenda, et massae materiali principium quoddam superius, et ut sic dicam formale addendum, quandoquidem omnes veritates rerum corporearum ex solis axiomatibus logisticis et geometricis, nempe de magno te parvo, toto et parte, figura et situ, colligi non possint, sed alia causa et effectu, actioneque et passione accedere debeant, quibus ordinis rerum rationes salventur. Id principium Formam, an ἐντελέχειαν, an Vim appellemus, non refert, modo meminerimus per solam virium notionem intelligibiliter explicari." (GM VI, 241-242). Cf. the quotations in footnote 65 and 67. For analogous statements see for instance what Leibniz's wrote some years before, in his *Discourses on Metaphysics* (1686): "Mais la force, or cause prochaine de ces changemens est quelque chose de plus reel, et il y a assez de fondement pour l'attribuer, à un corp plus qu'à l'autre; aussi n'est ce que par là qu'on peut connoistre à qui le mouvement appartient d'avantage. Or cette force est quelque chose de different de la grandeur, de la figure et du mouvement, et on peut juger par là que tout ce qui est conçu dans les corps ne consiste pas uniquement dans l'étendue et dans ses modifications, comme nos modernes se le persuadent. Ainsi nous sommes encor obligés de rétablir quelques estres ou formes, qu'il ont bannies. Et il paroist de plus en plus quoyque tous les phenomenes particuliers de la nature se puissent expliquer mathematiquement ou mechaniquement par ceux qui les entendent, qui neantmoins les

What has to be dealt with using mathematics is only the phenomenal "expression" of metaphysical forces, namely physical bodies and their interactions. As usual, dealing with physics, one must be aware that mathematical descriptions are idealisations, useful but necessarily approximate. For instance, Leibniz defines matter as a contiguous stuff but, as we mentioned in phase IV, contiguity is not a mathematical property. The better mathematical description of a contiguous physical body is still a continuous mathematical solid, where continuity still means, as in phase I, potentially infinite divisibility, and indefiniteness of parts.[105]

More generally, this holds for each physical notion involving extension – whether matter, time or motion – which therefore admit a traditional geometric representation, and can be dealt with using the traditional tools of Greek geometry, and in particular the method of exhaustion, without any need for infinitesimal quantities. Of course, since one can solve traditional geometrical problems using infinitesimal calculus, one can solve cinematic problems using $x$, $t$ and $dx$, $dt$. But while $x$ and $t$ have a direct counterpart in a measurable, and geometrically representable extensions, $dx$ and $dt$ have no direct physical counterparts. In other words, even if space and time are only potentially infinitely divisible (without ever attaining infinitesimals), we can fictively introduce infinitesimals of space, or of time.

As we argued at the end of phase IV, this is finally the sense in which we interpret Leibniz's idea of fictionality: a mathematical notion is a fiction when it lacks a traditional geometrical representation, and therefore has no physical counterpart.[106]

We conclude our analysis by observing that phase V offers a new perspective on the usefulness of fictions like infinitesimals, which Leibniz employs in an analogical way to speak of otherwise inexpressible metaphysical notions. For instance, in the *Specimen Dynamicus* (GM VI, 238), a *nisus* is described as an elementary, infinitely small solicitation, from the infinite repetition of which *impetus* arises. An analogous relation is described between the dead derivative force (*vis mortua*), to which the elementary solicitation belongs, and the living derivative force (*vis viva*).

Our analysis of Leibniz's approach to force is compatible with the analysis by Rutherford, who wrote:

> In the end, I propose, it is Leibniz's view that, in and of itself, force is not a mathematically representable property. To understand the sense in which force is real, we must turn to a different theoretical framework altogether, that of metaphysics, wherein force is represented as a modification of a substantial power, or principle of change. (Rutherford 2008 p. 257)

---

principes generaux de la nature corporelle et de la mecanique même sont plustot metaphysiques que Geometriques, et appartiennent plustost à quelques formes ou natures indivisibles comme causes des apparences qu'à la masse corporelle et étendue " (GP IV, 444. English translation in AG 51-52).

[105] "A continuum is a whole whose parts are indefinite, and so is the space if we abstract everything that is in it." [*Divisio terminorum*…1683 A6.4. 565]. The relation between continuity and indefiniteness of parts is at the basis of Aristotle's discussion of the potentially infinite division of the continuum, in *Physics* III.4-8 and VI.1, a discussion which Leibniz continued to appreciate far beyond his "aristotelian" phase I: "Aristote a fort bien expliqué le plein et la division du continu contre les atomistes" [*Remarques sur la doctrine Cartesienne* 1689. A6.4. 2047].

[106] In fact the introduction of infinitesimals, as well as other fictional mathematical notions, contributed to the release of algebra from geometry: "Leibniz avait adopté une nouvelle manière de pratiquer la géometrie. Le calcul entre ses mains n'était plus la traduction des figures et des problèmes mais devenait le vrai point de départ" (Parmentier 1995 p. 104). On this point see also Sherry-Katz 2014: "In this respect Leibniz recognized a hidden power of algebra, realizing that it was more than a short hand for what could be expressed geometrically. As impressed as Leibniz was by imaginary expression for real roots, he was more impressed by the power of algebra, and more generally, symbolic thinking, to unify disparate domains" (p.170). Even though he starts from different premises – that infinitesimals are compendia for Archimedean procedures and foreshadow the epsilontic approach to the notion of limit – an analogous conclusion is reached by Breger 2008 (pp. 195-198).

The connection with infinitesimal calculus together with the necessity of taking it as an analogy is acknowledged by Leibniz himself, and made explicit in a 1698 letter to De Volder:

> By analogy with geometry, or my analysis, solicitations are as $dx$, speeds as $x$, and [living] forces are as $xx$, or $\int x\,dx$.[107]

And this is in accordance with Leibniz's belief that metaphysical notions elude our human faculties, and the way in which we can speak of them necessarily involves reasoning by *analogy*:

> On ne doit point s'étonner que je tâche d'éclaircir ces choses par des *comparaisons* prises des mathématiques pures, où tout va dans l'ordre, et où il y a moyen de les démêler par une méditation exacte, qui nous fait jouir, pour ainsi dire, de la vue des idées [*Théodicée*, § 242].[108]
>
> ("One should not be surprised that I try to clarify these things by comparisons taken from pure mathematics, where everything goes in order, and where there is a way to disentangle them by a precise meditation, which makes us benefit, so to speak, from the view of ideas.")
> Je trouvay donc que leur nature consiste dans la force, et que de cela s'ensuit quelque chose d'*analogique* au sentiment et à l'appetit; et qu'ainsi il falloit les concevoir *à l'imitation* de la notion que nous avons des ames [*Systeme nouveau* (1695)].[109]
> ("I therefore found that their nature consists in the force, and that from this it follows something analogous with regard to feeling and appetite; and that thus it was necessary to conceive them in imitation of the notion that we have of souls.")

CONCLUSION

In this article, we have attempted to trace the evolution of Leibniz's thinking regarding infinitesimals, both in mathematics and in physics.

We have observed how Leibniz's increasing mathematical awareness led him to gradually abandon his initial scholastic position, centered around the idea of the potential infinite divisibility of the continuum, in favor of more informed reevaluations, that involve the notion of unassignability.

We have divided Leibniz's journey into five phases, which can be summarized as follows:
Phase I: both in mathematics and in physics Leibniz adopts a traditional syncategorematic approach: both extension and matter are continuous, which means that they are potentially infinitely divisible, and every obtained part is assignable. The problem arises of the actuality of the elementary constituents of physical matter.
Phase II: In an attempt to resolve the problem, Leibniz separates the treatment of physics from that of mathematics. He redefines matter as a continuous stuff infinitely divided – and not only divisible

---

[107] "Ut ita secundum analogiam Geometriae seu analysis nostrae solicitationes sint ut $dx$, celeritates ut $x$, vires ut $xx$ seu ut $\int x\,dx$" [GP II, 156].

[108] GP VI, 261-262. "je ne prends les mots dans ce sens extraordinaire que parceque je ne trouve point des termes plus propres pour m'exprimer" (*Additions à l'Explication du systeme nouveau touchant l'union de l'ame et du corps*, 1698. GP IV, 582). "L'Idée positive de cette substance simple ou Force primitive est toute trouvée, puisque elle doit toujours avoir en elle un progres reglé de perceptions, suivant l'Analogie qu'elle doit avoir avec notre ame" (Leibniz to Lady Masham 30 june 1704. GP III, 357). On the role of analogy in Leibniz's philosophy see Rescher 2013, ch. 7 and Pasini 2016.

[109] GP IV 479, translated in AG 139. Italics ours. Of course, it is not only the notion of force that cannot be translated in mathematical terms: think about the notion of simple substance, or that of hypercategorematic infinite, or even possible worlds etc... and the mathematical analogies which, this notwithstanding, Leibniz employs to illustrate them: "Vous savés donc que lorsque les conditions d'un point qu'on demande, ne le determinent pas assés, et qu'il y en a une infinité, ils tombent tous dans ce que les Geometres appellent un lieu, et ce lieu au moins (qui est souvent une Ligne) sera determiné. Ainsi vous pouvés vous figurer une suite reglée de Mondes, qui contiendront tous et seuls le cas dont il s'agit, et en varieront les circonstances et les consequences. Mais si vous posés un cas qui ne differe du monde actuel que dans une seule chose definie et dans ses suites, un certain monde determiné vous repondra" (*Théodicée*, § 414, GP VI, 362-363).

– into individual parts. Both actuality and unity of these parts are due to their conatus, or unextended amounts of motion. The geometrical counterpart of this physical situation is a hybrid in which the infinite divisibility of the continuum has to coexist with the actuality of unextended "fat points", which Leibniz identifies with Cavalieri's indivisibles. In this phase unassignability is accepted, without any further elaboration.

Phase III: Perhaps dissatisfied with the conclusions reached in phase II, Leibniz introduces some variations. The most significant one is the inversion of perspective in physics: discreteness is considered a foundational notion in relation to continuity: physical bodies do not come from an infinitely divided continuous stuff, but from the aggregation of discrete basic elements, that as in phase II are "glued" together by their conatus of motion. In mathematics the situation is similar, with the only difference being that now points are extended, although unassignable. Apparently, in this phase to keep together extended but unassignable quantities and assignable ones is not a problem for Leibniz. Unassignable quantities are not called fictions, nor referred to using cognate verbal forms.

Phase IV: Leibniz delves deeper into the concept of unassignability and further separates mathematical inquiry from physical inquiry. Unassignable quantities are useful in mathematics but cannot have a physical counterpart. In fact, he rejects the possibility that physical quantities are not all comparable to each other: a world containing some portion of matter which has an unassignable ratio to another portion would be less perfect than a world where all portions are comparable. Leibniz starts referring to unassignable quantities as "fictions" and links them to other "imaginary" notions, such as imaginary numbers, that are non-contradictory but do not admit any traditional geometrical representation, nor a physical counterpart. In physics, metaphysical explanations reappear alongside the mathematical ones, and in some cases, as a replacement for them.

Phase V: The separation between mathematics, physics and metaphysics is gradually brought to completion. In particular, mathematics ceases to be the foundation of physics, and geometry ceases to be its natural description: the true essence of physics is metaphysics, and as such it is non representable. In particular, being extraneous to the category of quantity, metaphysical objects can be thought of as unassignable in a stronger sense, and without problematic consequences.

We are far away from the reassuring "syncategorematic" vision of the beginning: mathematical awareness had led Leibniz to remarkable achievements in mathematics, but it had also compelled him to introduce greater complexity in physics.